\documentclass[reqno,12pt,a4paper]{amsart}

\usepackage{mathrsfs}
\usepackage{amssymb}
\usepackage{amsfonts}
\usepackage{latexsym}
\usepackage{amsthm}

\usepackage{enumerate}
\usepackage{titletoc}
\numberwithin{equation}{section}

\usepackage{graphicx}
\usepackage{tikz}

\begin{document}


\renewcommand{\theequation}{\arabic{section}.\arabic{equation}}
\theoremstyle{plain}
\newtheorem{theorem}{\bf Theorem}[section]
\newtheorem{lemma}[theorem]{\bf Lemma}
\newtheorem{corollary}[theorem]{\bf Corollary}
\newtheorem{proposition}[theorem]{\bf Proposition}
\newtheorem{definition}[theorem]{\bf Definition}
\newtheorem*{definition*}{\bf Definition}
\newtheorem*{example}{\bf Example}
\theoremstyle{remark}
\newtheorem*{remark}{\bf Remark}

\def\a{\alpha}  \def\cA{{\mathcal A}}     \def\bA{{\bf A}}  \def\mA{{\mathscr A}}
\def\b{\beta}   \def\cB{{\mathcal B}}     \def\bB{{\bf B}}  \def\mB{{\mathscr B}}
\def\g{\gamma}  \def\cC{{\mathcal C}}     \def\bC{{\bf C}}  \def\mC{{\mathscr C}}
\def\G{\Gamma}  \def\cD{{\mathcal D}}     \def\bD{{\bf D}}  \def\mD{{\mathscr D}}
\def\d{\delta}  \def\cE{{\mathcal E}}     \def\bE{{\bf E}}  \def\mE{{\mathscr E}}
\def\D{\Delta}  \def\cF{{\mathcal F}}     \def\bF{{\bf F}}  \def\mF{{\mathscr F}}
\def\c{\chi}    \def\cG{{\mathcal G}}     \def\bG{{\bf G}}  \def\mG{{\mathscr G}}
\def\z{\zeta}   \def\cH{{\mathcal H}}     \def\bH{{\bf H}}  \def\mH{{\mathscr H}}
\def\e{\eta}    \def\cI{{\mathcal I}}     \def\bI{{\bf I}}  \def\mI{{\mathscr I}}
\def\p{\psi}    \def\cJ{{\mathcal J}}     \def\bJ{{\bf J}}  \def\mJ{{\mathscr J}}
\def\vT{\Theta} \def\cK{{\mathcal K}}     \def\bK{{\bf K}}  \def\mK{{\mathscr K}}
\def\k{\kappa}  \def\cL{{\mathcal L}}     \def\bL{{\bf L}}  \def\mL{{\mathscr L}}
\def\l{\lambda} \def\cM{{\mathcal M}}     \def\bM{{\bf M}}  \def\mM{{\mathscr M}}
\def\L{\Lambda} \def\cN{{\mathcal N}}     \def\bN{{\bf N}}  \def\mN{{\mathscr N}}
\def\m{\mu}     \def\cO{{\mathcal O}}     \def\bO{{\bf O}}  \def\mO{{\mathscr O}}
\def\n{\nu}     \def\cP{{\mathcal P}}     \def\bP{{\bf P}}  \def\mP{{\mathscr P}}
\def\r{\varrho} \def\cQ{{\mathcal Q}}     \def\bQ{{\bf Q}}  \def\mQ{{\mathscr Q}}
\def\s{\sigma}  \def\cR{{\mathcal R}}     \def\bR{{\bf R}}  \def\mR{{\mathscr R}}
\def\S{\Sigma}  \def\cS{{\mathcal S}}     \def\bS{{\bf S}}  \def\mS{{\mathscr S}}
\def\t{\tau}    \def\cT{{\mathcal T}}     \def\bT{{\bf T}}  \def\mT{{\mathscr T}}
\def\f{\phi}    \def\cU{{\mathcal U}}     \def\bU{{\bf U}}  \def\mU{{\mathscr U}}
\def\F{\Phi}    \def\cV{{\mathcal V}}     \def\bV{{\bf V}}  \def\mV{{\mathscr V}}
\def\P{\Psi}    \def\cW{{\mathcal W}}     \def\bW{{\bf W}}  \def\mW{{\mathscr W}}
\def\o{\omega}  \def\cX{{\mathcal X}}     \def\bX{{\bf X}}  \def\mX{{\mathscr X}}
\def\x{\xi}     \def\cY{{\mathcal Y}}     \def\bY{{\bf Y}}  \def\mY{{\mathscr Y}}
\def\X{\Xi}     \def\cZ{{\mathcal Z}}     \def\bZ{{\bf Z}}  \def\mZ{{\mathscr Z}}
\def\O{\Omega}

\newcommand{\mc}{\mathscr {c}}

\newcommand{\gA}{\mathfrak{A}}          \newcommand{\ga}{\mathfrak{a}}
\newcommand{\gB}{\mathfrak{B}}          \newcommand{\gb}{\mathfrak{b}}
\newcommand{\gC}{\mathfrak{C}}          \newcommand{\gc}{\mathfrak{c}}
\newcommand{\gD}{\mathfrak{D}}          \newcommand{\gd}{\mathfrak{d}}
\newcommand{\gE}{\mathfrak{E}}
\newcommand{\gF}{\mathfrak{F}}           \newcommand{\gf}{\mathfrak{f}}
\newcommand{\gG}{\mathfrak{G}}           
\newcommand{\gH}{\mathfrak{H}}           \newcommand{\gh}{\mathfrak{h}}
\newcommand{\gI}{\mathfrak{I}}           \newcommand{\gi}{\mathfrak{i}}
\newcommand{\gJ}{\mathfrak{J}}           \newcommand{\gj}{\mathfrak{j}}
\newcommand{\gK}{\mathfrak{K}}            \newcommand{\gk}{\mathfrak{k}}
\newcommand{\gL}{\mathfrak{L}}            \newcommand{\gl}{\mathfrak{l}}
\newcommand{\gM}{\mathfrak{M}}            \newcommand{\gm}{\mathfrak{m}}
\newcommand{\gN}{\mathfrak{N}}            \newcommand{\gn}{\mathfrak{n}}
\newcommand{\gO}{\mathfrak{O}}
\newcommand{\gP}{\mathfrak{P}}             \newcommand{\gp}{\mathfrak{p}}
\newcommand{\gQ}{\mathfrak{Q}}             \newcommand{\gq}{\mathfrak{q}}
\newcommand{\gR}{\mathfrak{R}}             \newcommand{\gr}{\mathfrak{r}}
\newcommand{\gS}{\mathfrak{S}}              \newcommand{\gs}{\mathfrak{s}}
\newcommand{\gT}{\mathfrak{T}}             \newcommand{\gt}{\mathfrak{t}}
\newcommand{\gU}{\mathfrak{U}}             \newcommand{\gu}{\mathfrak{u}}
\newcommand{\gV}{\mathfrak{V}}             \newcommand{\gv}{\mathfrak{v}}
\newcommand{\gW}{\mathfrak{W}}             \newcommand{\gw}{\mathfrak{w}}
\newcommand{\gX}{\mathfrak{X}}               \newcommand{\gx}{\mathfrak{x}}
\newcommand{\gY}{\mathfrak{Y}}              \newcommand{\gy}{\mathfrak{y}}
\newcommand{\gZ}{\mathfrak{Z}}             \newcommand{\gz}{\mathfrak{z}}

\def\ve{\varepsilon}   \def\vt{\vartheta}    \def\vp{\varphi}    \def\vk{\varkappa}

\def\A{{\mathbb A}} \def\B{{\mathbb B}} \def\C{{\mathbb C}}
\def\dD{{\mathbb D}} \def\E{{\mathbb E}} \def\dF{{\mathbb F}} \def\dG{{\mathbb G}} \def\H{{\mathbb H}}\def\I{{\mathbb I}} \def\J{{\mathbb J}} \def\K{{\mathbb K}} \def\dL{{\mathbb L}}\def\M{{\mathbb M}} \def\N{{\mathbb N}} \def\O{{\mathbb O}} \def\dP{{\mathbb P}} \def\R{{\mathbb R}} \def\dQ{{\mathbb Q}}
\def\S{{\mathbb S}} \def\T{{\mathbb T}} \def\U{{\mathbb U}} \def\V{{\mathbb V}}\def\W{{\mathbb W}} \def\X{{\mathbb X}} \def\Y{{\mathbb Y}} \def\Z{{\mathbb Z}}

\newcommand{\1}{\mathbbm 1}
\newcommand{\dd}    {\, \mathrm d}



\def\la{\leftarrow}              \def\ra{\rightarrow}            \def\Ra{\Rightarrow}
\def\ua{\uparrow}                \def\da{\downarrow}
\def\lra{\leftrightarrow}        \def\Lra{\Leftrightarrow}


\def\lt{\biggl}                  \def\rt{\biggr}
\def\ol{\overline}               \def\wt{\widetilde}
\def\no{\noindent}


\let\ge\geqslant                 \let\le\leqslant
\def\lan{\langle}                \def\ran{\rangle}
\def\/{\over}                    \def\iy{\infty}
\def\sm{\setminus}               \def\es{\emptyset}
\def\ss{\subset}                 \def\ts{\times}
\def\pa{\partial}                \def\os{\oplus}
\def\om{\ominus}                 \def\ev{\equiv}
\def\iint{\int\!\!\!\int}        \def\iintt{\mathop{\int\!\!\int\!\!\dots\!\!\int}\limits}
\def\el2{\ell^{\,2}}             \def\1{1\!\!1}
\def\sh{\sharp}
\def\wh{\widehat}

\def\all{\mathop{\mathrm{all}}\nolimits}
\def\where{\mathop{\mathrm{where}}\nolimits}
\def\as{\mathop{\mathrm{as}}\nolimits}
\def\Area{\mathop{\mathrm{Area}}\nolimits}
\def\arg{\mathop{\mathrm{arg}}\nolimits}
\def\adj{\mathop{\mathrm{adj}}\nolimits}
\def\const{\mathop{\mathrm{const}}\nolimits}
\def\det{\mathop{\mathrm{det}}\nolimits}
\def\diag{\mathop{\mathrm{diag}}\nolimits}
\def\diam{\mathop{\mathrm{diam}}\nolimits}
\def\dim{\mathop{\mathrm{dim}}\nolimits}
\def\dist{\mathop{\mathrm{dist}}\nolimits}
\def\Im{\mathop{\mathrm{Im}}\nolimits}
\def\Iso{\mathop{\mathrm{Iso}}\nolimits}
\def\Ker{\mathop{\mathrm{Ker}}\nolimits}
\def\Lip{\mathop{\mathrm{Lip}}\nolimits}
\def\rank{\mathop{\mathrm{rank}}\limits}
\def\Ran{\mathop{\mathrm{Ran}}\nolimits}
\def\Re{\mathop{\mathrm{Re}}\nolimits}
\def\Res{\mathop{\mathrm{Res}}\nolimits}
\def\res{\mathop{\mathrm{res}}\limits}
\def\sign{\mathop{\mathrm{sign}}\nolimits}
\def\supp{\mathop{\mathrm{supp}}\nolimits}
\def\Tr{\mathop{\mathrm{Tr}}\nolimits}
\def\AC{\mathop{\rm AC}\nolimits}
\def\BBox{\hspace{1mm}\vrule height6pt width5.5pt depth0pt \hspace{6pt}}


\newcommand\nh[2]{\widehat{#1}\vphantom{#1}^{(#2)}}
\def\dia{\diamond}

\def\Oplus{\bigoplus\nolimits}




\def\qqq{\qquad}
\def\qq{\quad}
\let\ge\geqslant
\let\le\leqslant
\let\geq\geqslant
\let\leq\leqslant

\newcommand{\ca}{\begin{cases}}
\newcommand{\ac}{\end{cases}}
\newcommand{\ma}{\begin{pmatrix}}
\newcommand{\am}{\end{pmatrix}}
\renewcommand{\[}{\begin{equation}}
\renewcommand{\]}{\end{equation}}
\def\bu{\bullet}

\title[{Periodic Dirac operator on the half-line}]
{Periodic Dirac operator on the half-line}

\date{\today}

\author[Evgeny Korotyaev]{Evgeny Korotyaev}
\address{Department of Mathematical Analysis,  Saint-Petersburg State University,   Universitetskaya nab. 7/9, St.
Petersburg, 199034, Russia, \ korotyaev@gmail.com, \ e.korotyaev@spbu.ru}
\author[Dmitrii Mokeev]{Dmitrii Mokeev}
\address{Saint-Petersburg State University,   Universitetskaya nab. 7/9, St.
Petersburg, 199034, Russia, \ mokeev.ds@yandex.ru}


\begin{abstract}
We consider the Dirac operator with a periodic potential on the
half-line with the Dirichlet boundary condition at zero. Its
spectrum consists of an absolutely continuous part plus at most one
eigenvalue in each open gap. The Dirac resolvent admits a
meromorphic continuation onto a two-sheeted Riemann surface with a
unique simple pole on each open gap: on the first sheet (an
eigenvalue) or on the second sheet (a resonance). These poles are
called states and there are no other poles. If the potential is
shifted by real parameter $t$, then the continuous spectrum does not
change but the states can change theirs positions. We prove that
each state is smooth and in general, non-monotonic function of $t$.
We prove that a state is a strictly monotone function of $t$ for a
specific potential. Using these results we obtain formulas to
recover potentials of special forms.

\end{abstract}

\maketitle


{\it \footnotesize Dedicated to the memory of Boris Dubrovin
(1950-2019)}

\section{Introduction} \label{p0}

We consider the Dirac operator $H^{\pm}_t$, $t \in \R$, acting on
$L^2 (\R_{\pm}, \C^2)$ and given by
$$
        H^{\pm}_t f = J f'+ V_t f,\qq f(x) = \ma f_1(x) \\ f_2(x) \am,
        \qq x \in \R_\pm,\qq J = \ma 0 & 1 \\ -1 & 0 \am,
$$
under the Dirichlet boundary condition $f_1(0)=0$, where $t \in \R$
is a shift. Here the $2\ts 2$ matrice $V_t(x)=V(x+t)$ is 1-periodic
and belongs to the Hilbert space $\cP$ defined  by
$$
\cP = \rt\{ \, \left. V= \ma q_1 & q_2 \\ q_2 & -q_1 \am \, \right| \, q_1,q_2
\in L_{real}^2(\T) \, \rt\},\qq         \T = \R/\Z,
$$
equipped with the norm $\left\| V \right\|_{\cP}^2 ={1\/2} \int_{\T}
 \Tr V^2(x)dx$.  The operator $H^{\pm}_t$ is self-adjoint, its spectrum
consists of an absolutely continuous part $\s_{ac}(H^{\pm}_t) =
\s_{ac}(H^{\pm}_0) = \cup_{n \in \Z} \s_n$ plus  at most one eigenvalue in
each non-empty gap $\g_n$, where the bands $\s_n$ and gaps $\g_n$
are given by
$$
\s_n = [\a_{n-1}^+,\a_n^-],\qq \g_n = (\a_{n}^-,\a_n^+),\qq \text{satisfying}\qq
\a_{n-1}^+ < \a_n^- \leq \a_n^+\qq \forall \, n \in \Z.
$$
 The sequence $\a_n^\pm$, $n \in \Z$,
is the spectrum of the equation $J y' + Vy = \l y$ with the
condition of 2-periodicity, $y(x+2)=y(x)$ $(x\in \R)$. If some  gap
degenerates, $\g_n= \es $, then the corresponding bands $\s_{n} $ and
$\s_{n+1}$ touch. This happens when $\a_n^-=\a_n^+$; this number is
then a double eigenvalue of the 2-periodic problem.  Generally, the
eigenfunctions corresponding to eigenvalues $\a_{2n}^{\pm}$ are
1-periodic, those for $\a_{2n+1}^{\pm}$ are 1-anti-periodic in the
sense that $y(x+1)=-y(x)$ $(x\in\R)$.

For the Dirac operator we  introduce the two-sheeted Riemann surface
$\L$ obtained by joining the upper and lower rims of two copies of
the cut plane $\C\sm\s_{ac}(H^+_0)$ in the usual (crosswise) way (see Fig.~1 below).
We denote the $n$-th gap on the first, physical, sheet $\L_1$ by
$\g_n^{(1)}$ and its counterpart on the second, nonphysical, sheet
$\L_2$ by $\g_n^{(2)}$, and set a circle gap $\g_n^c
=\ol \g_n^{(1)} \cup \ol \g_n^{(2)}$.


Let $f_{\pm}(\l,t)=((H^{\pm}_t-\l)^{-1}\eta,\eta)$, $(\l,t) \in \L_1
\ts \R$, for  some $\eta \in C_o^\iy(\R_\pm)$, $\eta \ne 0$. For
each $t \in \R$ the function $f_{\pm}(\cdot,t)$ admits a meromorphic
continuation from $\L_1$ onto the Riemann surface $\L$ with at most
one simple pole on each circle gap $\g_n^c$, its position depends on
$t$ and there are no other poles.
 A pole on the first sheet is an eigenvalue of $H^{\pm}_t$. A pole on
the second sheet is called a {\it resonance}. A point $\l_o=\a_n^\pm$, $n \in \Z$, is called a {\it virtual state} if the function $z
\mapsto f(\l_o+z^2)$ has a pole at $0$. A point $\l_o \in \L$ is
called a {\it state} if it is either an eigenvalue or a resonance or
a virtual state.  It is important that the projection of the pole on
the complex plane coincides with the Dirichlet eigenvalue $\m_n(t)$ of the
equation $J y' + V_t y = \l y$ on the unit interval with the boundary conditions
$y_1(0)=y_1(1)=0$, where $y = \ma y_1 & y_2 \am^{\top}$.


We consider the corresponding problem for a Schr{\"o}dinger
operator on $L^2(\R_{\pm})$ given by
$$
h^{\pm}_t y=-y'' + q(\cdot+t)y,\qq y(0) = 0,\qq q \in L^2(\T).
$$
The operator $h^{\pm}_t$ is self-adjoint, its spectrum consists of
an absolutely continuous part $\s_{ac}(h^{\pm}_t) = \s_{ac}(h^{+}_0)
= \cup_{n \ge 1} s_n$ plus  at most one eigenvalue in each non-empty
gap $g_n$, where the bands $s_n$ and gaps $g_n$ are defined  by
$$
s_n = [a_{n-1}^+,a_n^-],\qq g_n = (a_{n}^-,a_n^+),\qq  {\rm
satisfying}\qq a_{n-1}^+ < a_n^- \leq a_n^+\qq \forall \, n \ge 1.
$$

The sequence $a_0^+, a_n^\pm$, $n\ge 1$,
is the spectrum of the equation $-y'' + qy = \l y$ with the
condition of 2-periodicity, $y(x+2)=y(x)$ $(x\in \R)$. If some  gap
degenerates, $g_n=\es $, then the corresponding bands $s_{n} $ and
$s_{n+1}$ touch. This happens when $a_n^-=a_n^+$; this number is
then a double eigenvalue of the 2-periodic problem.  Generally, the
eigenfunctions corresponding to eigenvalues $a_{2n}^{\pm}$ are
1-periodic, those for $a_{2n+1}^{\pm}$ are 1-anti-periodic in the
sense that $y(x+1)=-y(x)$ $(x\in\R)$.

For the operator $h^{\pm}_t$
 we  introduce the two-sheeted Riemann surface
$\cL$ obtained by joining the upper and lower rims of two copies of
the cut plane $\C\sm\s_{ac}(h^{+}_0)$ in the usual (crosswise) way.
We denote the $n$-th gap on the first, physical, sheet $\cL_1$ by
$g_n^{(1)}$ and its counterpart on the second, nonphysical, sheet
$\cL_2$ by $g_n^{(2)}$, and set a circle gap $g_n^c
=\ol g_n^{(1)}\cup \ol g_n^{(2)}$.

For the Schr{\"o}dinger operator one can define states
 (eigenvalues, resonances and  virtual states) as we have
defined them above for the Dirac operator. Thus in each open gap
$g_n^c$ there exists a unique state $\r_n(t)$ of $h^{+}_t$ and its projections on the complex plane coincide with the
Dirichlet eigenvalues on the unit interval, i.e. $-y'' +
q(x+t)y=\r_n(t) y,\ y(0) =y(1)=0$ (see \cite{KS12, Zh}). We introduce
the Sobolev spaces $\cH^\a(\T)$, $\a \geq 0$, of all functions  $f,
f^{(\a)} \in
L^2(\T)$. The states $\r_n(t)$ satisfy:

\no $\bu$  {\it $\r_n(t)$ is an eigenvalue of $h^{\pm}_t$
    iff
$\r_n(t)_*$ is a resonance of $h^{\mp}_t$, where $\r_n(t)_* \in \cL$ is projection of $\r_n(t)$ to the other sheet.

\no    $\bu$ $\r_n(\cdot) \in \cH^2(\T)$ and  $\r_n(t)$ runs
strictly monotonically around $g_n^c$, i.e. it moves in one
direction along $g_n^{(1)}$ and in the opposite direction along
$g_n^{(2)}$, changing sheets when it hits $a_n^+$ or $a_n^-$ and
making $n$ complete revolutions around $g_n^c$ when $t$ runs through
$[0,1)$.

\no    $\bu$ $\r_n(t)$ is a solution of the Dubrovin equation
$$
\r_n'(t) = \frac{n^2 \pi^2 \sqrt{D^2(\r_n(t)) - 1}}{\prod_{j \neq n}
\left( \frac{\r_j(t)-\r_n(t)}{j^2 \pi^2} \right)},\ \ n
=1,2,3,\ldots,
$$
        where $(-1)^{j} \sqrt{D^2(\l) - 1} < 0$ for $\l \in g_n^{(j)}$, $j=1,2$, and
        the following factorization holds true:}
$$
D^2(\l)-1 = (a_0 - \l)\prod_{j \geq 1} \frac{(a_j^- - \l)(a_j^+ -
\l)}{j^4 \pi^4}.
$$

In order to study the KdV equation on the circle Dubrovin \cite{D}
considered the dynamics of $\r_n(t)$ for  finite band potentials.
This problem in connection with the inverse problem for the periodic
Schr{\"o}dinger operator was considered by Trubowitz \cite{trub77}
in the case of smooth potential (see also \cite{Lev}) and by
Korotyaev \cite{Kor99} in the case of potential from $L^2(\T)$ (see also \cite{KS12,
Zh}). These results are used to construct the gap-length mapping for
inverse problem for the operator $h^{+}_0$ in \cite{Kor99}.

Such operators are used to study defects in a crystal. For example, in \cite{Kor01a}
this operator was used to investigate a one-dimensional
Schr{\"o}dinger operator with a dislocation in periodic media (see
also \cite{DPR09, D19, Fef, Hemp, HK11b, HKS15}). Using Schr{\"o}dinger operator on the
half-line, Korotyaev \cite{Kor05} studied junction of two
one-dimesional periodic potentials for Schr{\"o}dinger operator, and
in particular junction of crystal and vacuum.

Our main goal is to obtain similar results for the Dirac operator
$H^{+}_t$. We will show that in each circle open gap there exists a
unique state of $H^{+}_t$. Each state is a function from $\cH^1(\T)$,
and its projection on the plane coincides with the Dirichlet
eigenvalues on the unit interval. We show that the states can be
non-monotonic  or monotonic functions of $t$ under some restrictions
on potentials. Moreover, we will obtain local asymptotic of states
and will get some results about inverse problem using this
asymptotic.

The periodic Dirac operator on the half-line has been studied in many papers (see literature in \cite{BES13}). The Dirichlet eigenvalues for the Dirac operator was considered
as function of the shift parameter in connection with the nonlinear
Schr{\"o}dinger equation. In the finite band case, the motion of the
Dirichlet eigenvalues was considered in \cite{Its76, KotIts14,
Prev85}. The Dubrovin equation with a  smooth potential was studied
in \cite{GreGui}, \cite{BGGK}.  It was shown that the image of each
Dirichlet eigenvalue $\m_n(t)$ covers the gap $\g_n$ when $t$ runs
through $[0,1]$ \cite{BGGK}. The inverse problem for the periodic
Dirac operator in terms of gaps lengths was solved in the paper
\cite{Kor05b}. This inverse problem in terms of the
Marchenko-Ostrovki mapping was solved in the papers \cite{Kor01,
Mi}.

To describe the states $\r_n(t)$ of the Schr{\"o}dinger operator as
functions of the shift in the paper \cite{Kor99} the implicit
function theorem is used. We cannot use this version of the implicit
function theorem since solutions of the Dirac equation are less
smooth than solutions of the Schr{\"o}dinger equation. Thus, we use
a specific version of the implicit function theorem.

This paper is organized as follows:\\
 In Section \ref{p1} we
formulate our main results.\\
 In Section \ref{p2} we consider the periodic Dirac operator on the half-line.
\\ In Section \ref{p3} we prove the main theorems.\\
In Section \ref{p6} (Appendix) we prove some technical lemmas, and the
specific version of the implicit function theorem.

\section{Main results} \label{p1}
\subsection{Definitions}
In order to formulate our main results, we shortly describe
properties of the Dirac operator on the half-line. We introduce the
fundamental vector-valued solutions $
        \vt(x,\l,t) = \ma \vt_1 \\ \vt_2 \am (x,\l,t)$,
        $\vp(x,\l,t) = \ma \vp_1 \\ \vp_2 \am (x,\l,t)
    $
of the shifted Dirac equation
\[ \label{p1e1}
        Jy'(x) + V(x+t) y(x) = \l y(x),\ \ (x,\, \l,\, t) \in \R \ts \C \ts \R,
\]
 satisfying the conditions
    $
        \vt(0,\l,t) = \ma 1 \\ 0 \am$, $\vp(0,\l,t) = \ma 0 \\ 1 \am
    $,
where  $V \in \cP$.

In Lemma \ref{p3l1} we show the standard fact that the resonances of
$H^{+}_t$ are the eigenvalues of $H_t^-$. We also prove that each
state of $H^{+}_t$ on $\L$ is over the Dirichlet eigenvalue
$\m_n(t)$ of the corresponding Dirichlet problem on the unit
interval. Therefore, following \cite{trub77}, we denote a state of
$H^{+}_t$ on $\L$ by $\m_n(t)$ as well as the Dirichlet eigenvalue.
It is also follows that each state of $H^{\pm}_t$ belongs to
$\g_n^c$. We show that if $\g_n^c$ is empty for some $n \in \Z$,
then $\m_n(t) = \a_n^{\pm}$ is not a state of $H^{\pm}_t$ for each
$t \in \T$. Moreover, we will show that for $\g_n^c \neq \es$
\begin{enumerate}[1)]
\item if $\left| \vp_2(1,\m_n(t),t) \right| < 1$, then $\m_n(t) \in \g_n^{(1)}$,
i.e. $\m_n(t)$ is an eigenvalue of $H^{+}_t$;
\item if $\left| \vp_2(1,\m_n(t),t) \right| > 1$, then $\m_n(t) \in \g_n^{(2)}$, i.e. $\m_n(t)$ is a resonance of $H^{+}_t$;
\item if $\left| \vp_2(1,\m_n(t),t) \right| = 1$, then $\m_n(t) = \a_n^+$ or $\m_n(t) = \a_n^-$, i.e. $\m_n(t)$ is a virtual state of $H^{+}_t$.
\end{enumerate}
\begin{remark}  These quantities are used to construct  the gap length
mapping in \cite{Kor05b}.
\end{remark}
    \begin{figure}[h]
        \begin{center}
            \begin{tikzpicture}
            \usepgflibrary{decorations.markings}
            \usetikzlibrary{intersections}
            \usetikzlibrary{decorations}
            \usetikzlibrary{decorations.pathmorphing}
            \usetikzlibrary{decorations.markings}
            \usetikzlibrary{decorations.pathreplacing}
            \usetikzlibrary{decorations.shapes}
            \usetikzlibrary{topaths}
            \usetikzlibrary{calc}
            \usetikzlibrary{math}

            \tikzmath{
                real \height, \width, \ratio, \gap, \zone, \amp, \xs, \ys;
                coordinate \LD, \GLD, \a, \b, \c, \d, \e, \f, \g, \h, \i, \j, \k, \l, \m, \n;
                \height = 1.3; 
                \width = 12; 
                \ratio = 1.2; 
                \amp = 0.08; 
                \xs = 1.5; 
                \ys = 0; 
                \ss = -60; 
                \step = 4pt; 
                \LD = (0, 0); 
                \gap = \width / (2.6*\ratio + 3); 
                \zone = \ratio * \gap; 
                \GLD = (\LDx, \LDy) + (\width/2 - 2*\zone -1.5*\gap, 0); 
                \a = (\GLDx,\GLDy) + (\zone, 0);
                \b = (\LDx,\LDy) + (0,\amp);
                \c = (\LDx,\LDy) + (0, \height);
                \j = (\LDx,\LDy) + (0,-\amp);
                \i = (\LDx,\LDy) + (0, -\height);
                \d = (\cx,\cy) + (\width, 0);
                \e = (\bx,\by) + (\width, 0);
                \f = (\ax,\ay) + (3*\gap, 0) + (2*\zone, 0);
                \g = (\jx,\jy) + (\width, 0);
                \h = (\ix,\iy) + (\width, 0);
                \k = (\ax,\ay) + (\gap, 0);
                \l = (\kx,\ky) + (\zone, 0);
                \m = (\lx,\ly) + (\gap, 0);
                \n = (\mx,\my) + (\zone, 0);
            }

            \tikzset{
                border lines/.style={black, semithick},
                gap lines/.style={black, very thin},
                zone lines/.style={bend right = 10,looseness=0.8},
                filling/.style={fill=white, fill opacity=.6},
                first decor/.style={
                    postaction={
                        decoration={
                            crosses, raise=-3pt, shape size=3pt, pre length=0pt, shape sep = 7pt
                        },
                        decorate, draw
                    }},
                second decor/.style={
                    postaction={
                        decoration={
                            triangles, raise=-3pt, shape size=3pt, pre length=0pt, shape sep = 7pt
                        },
                        decorate, draw
                    }},
            }

            \begin{scope}[yshift=\ss, yslant=\ys, xslant=\xs]

                \draw[border lines, filling] (\a) to[zone lines] (\b) -- (\c) -- (\d) -- (\e)
                to[zone lines] (\f) to[zone lines] (\g) -- (\h) -- (\i) -- (\j) to[zone lines] cycle;

                \draw[gap lines] (\a) -- (\k) node[midway,below,black]{$\gamma_{n-1}^{(2)}$} to[move to] (\l) -- (\m) node[midway,below,black]{$\gamma_{n}^{(2)}$} node[pos=0.3,above,black]{$\mu_n(t)$} node[pos=0.3,fill=black, circle,minimum size=3pt,inner sep=0pt]{} to[move to] (\n) -- (\f) node[midway,black,below]{$\gamma_{n+1}^{(2)}$};

                \draw[border lines] (\l) to[zone lines] (\k) to[zone lines] cycle;
                \draw[border lines] (\n) to[zone lines] (\m) to[zone lines] cycle;

                \draw (\d) node[black, below = 11pt, left = 26pt]{$\Lambda_2$};

                \path[second decor] (\e) to[zone lines] (\f) to[move to] (\n) to[zone lines] (\m) to[move to] (\l) to[zone lines] (\k) to[move to] (\a) to[zone lines] (\b);
                \path[first decor] (\j) to[zone lines] (\a) to[move to] (\k) to[zone lines] node[pos = 0.5,black, below = 5pt, sloped, fill = white]{$\sigma_{n}$} (\l) to[move to] (\m) to[zone lines] node[pos = 0.5,black, below = 5pt, sloped, fill = white]{$\sigma_{n+1}$} (\n) to[move to] (\f) to[zone lines] (\g);
            \end{scope}

            \begin{scope}[yshift=0, yslant=\ys, xslant=\xs]

                \draw[border lines, filling] (\a) to[zone lines] (\b) -- (\c) -- (\d) -- (\e)
                to[zone lines] (\f) to[zone lines] (\g) -- (\h) -- (\i) -- (\j) to[zone lines] cycle;

                \draw[gap lines] (\a) -- (\k) node[midway,above,black]{$\gamma_{n-1}^{(1)}$} node[pos=0.7,below,black]{$\mu_{n-1}(t)$} node[pos=0.7,fill=black, circle,minimum size=3pt,inner sep=0pt]{} to[move to] (\l) -- (\m) node[midway,above,black]{$\gamma_{n}^{(1)}$} to[move to] (\n) -- (\f) node[midway,above,black]{$\gamma_{n+1}^{(1)}$} node[pos=0.4,below,black]{$\mu_{n+1}(t)$} node[pos=0.4,fill=black, circle,minimum size=3pt,inner sep=0pt]{};

                \draw[border lines] (\l) to[zone lines] (\k) to[zone lines] cycle;
                \draw[border lines] (\n) to[zone lines] (\m) to[zone lines] cycle;

                \draw (\d) node[black, below = 11pt, left = 26pt]{$\Lambda_1$};

                \path[first decor] (\e) to[zone lines] (\f) to[move to] (\n) to[zone lines] node[pos = 0.5,black, above = 5pt, sloped, fill = white]{$\sigma_{n+1}$} (\m) to[move to] (\l) to[zone lines] node[pos = 0.5,black, above = 5pt, sloped, fill = white]{$\sigma_{n}$} (\k) to[move to] (\a) to[zone lines] (\b);
                \path[second decor] (\j) to[zone lines] (\a) to[move to] (\k) to[zone lines] (\l) to[move to] (\m) to[zone lines] (\n) to[move to] (\f) to[zone lines] (\g);
            \end{scope}
        \end{tikzpicture}
        \label{fig1}
        \caption{Structure of the Riemann surface $\L$. Crosses and triangles show how rims on $\L_1$ and $\L_2$ are joined together. The dots denote positions of $\m_n(t)$.}
        \end{center}
    \end{figure}
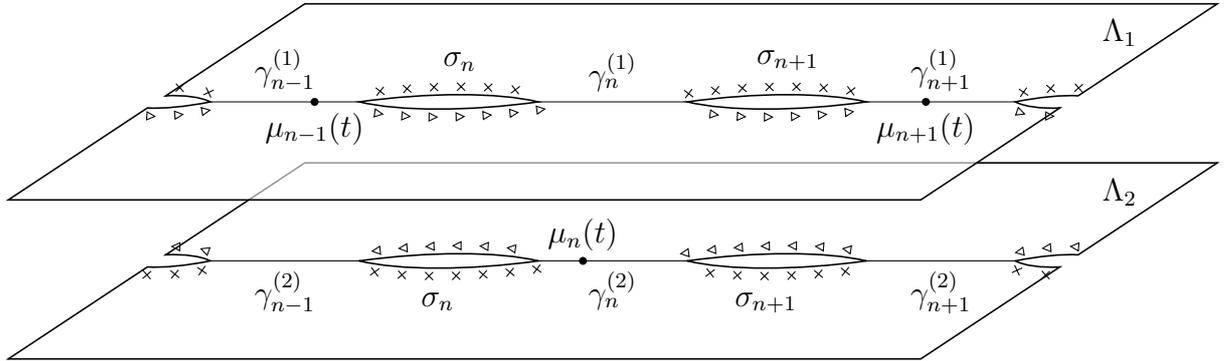

We also introduce the following functions
\[
\label{p1e11}
\begin{split}
\Omega(\l, t) = \frac{1 -
\vp_2(1,\l,t)^2}{\left\|\vp(\cdot,\l,t)\right\|^2},
        \quad \Omega_n^{\pm}(t) = \frac{2 M_n^{\pm} }{\left\|\vp(\cdot,\a_n^{\pm},t)\right\|^4},
        \quad \omega_n^{\pm}(t) = \frac{- 2 M_n^{\pm}}{\left\|\vp(\cdot,\a_n^{\pm},t)\right\|^2},
        \\
Q_n(t) = \int_{0}^{t} \left| q_1(\t) + \m_n^0 \right| d\t,\quad
U_n(t) = \left| t \right| + \left| t \right|^{1/2} \left| \int_{0}^t
\left| V(\t) + \m_n^0I_2 \right|^2 d\t \right|^{1/2},
    \end{split}
    \]
where $(\l,t,n) \in \C \times \R \times \Z$, $\m_n^0 = \m_n(0)$,
$I_2$ is the identity matrix $2 \times 2$, matrix norm $\left| \cdot
\right|$ is defined by $\left| A \right|^2 = \Tr A^* A$, norm of any
vector-valued function is defined by $\left\| \vp \right\|^2 =
\int_0^1 (\left|\vp_1(x)\right|^2 + \left|\vp_2(x)\right|^2)dx$,
effective masses $\pm M_n^{\pm} > 0$ are defined in Section \ref{p2}.
Here and below the dot denotes the derivative with respect to $t$,
i.e. $\dot{u} = du/dt$.

\subsection{Main theorems}

\begin{theorem} \label{p0t1}
Let $(t,V) \in \R \ts \cP$ and let a gap $\g_n$ be open for some $n \in
\Z$. Then there exists a unique state $\m_n(t)$ of the operator
$H^{+}_t$ in $\g_n^c$ and the following statements hold true:
\begin{enumerate}[i)]
\item $\m_n(\cdot) \in \cH^1(\T)$ and for almost all $t \in \T$ it satisfies the Dubrovin equation
\[ \label{p1e5}
\dot{\m}_n(t) = \left( q_1(t)+\m_n(t) \right) \Omega(\m_n(t), t).
\]
\item If $\m_n^0 \neq \a_n^{\pm}$, then the following asymptotic holds as $t \to 0$:
\[ \label{p1e3}
\m_n(t) = \m_n^0 + \Omega(\m_n^0, 0) \int_{0}^t (q_1(\t) + \m_n^0)
d\t + O\left(Q_n(t)U_n(t)\right).
            \]
\item If $\m_n^0 = \a_n^{\pm}$, then the following asymptotics hold as $t \to 0$:
\begin{align}
\m_n(t) = \m_n^0 + \Omega_n^{\pm}(0) \left( \int_{0}^t (q_1(\t) + \m_n^0) d\t
\right)^2 + O\left(Q_n(t)^2U_n(t)\right), \label{p1e4} \\
(-1)^n\vp_2(1,\m_n(t),t) = 1 + \omega_n^{\pm}(0) \int_{0}^t (q_1(\t)
+ \m_n^0) d\t + O\left(Q_n(t)U_n(t)\right). \label{p1e6}
\end{align}
\item Let $Y(\cdot,\a_n^{\pm},0)$ be $2$-periodic solution of equation (\ref{p1e1})
for $\l = \a_n^{\pm}$ and $t = 0$. Then $Y_1(t_0,\a_n^{\pm},0) = 0$ for some $t_0 \in [0,1)$
if and only if $\m_n(t_0) = \a_n^{\pm}$.
\end{enumerate}
    \end{theorem}
\begin{remark}
1) Recall that if $\g_n^c=\es$ for some $n \in \Z$, then
$\m_n(t) = \a_n^{\pm}$ is not a state of $H^{+}_t$ for each $t \in
\T$. We formulate asymptotics (\ref{p1e3}--6) as $t \to 0$. Shifting
$V$ by any $t_0 \in \R$, we determine similar asymptotics as $t \to
t_0$ (see (\ref{p4e41}--24)).

2) In the proof of Theorem \ref{p0t1} we used arguments from the
papers \cite{Kor99, trub77} and the specific implicit function
theorem. In the case of the Dirac operator $\m_n(\cdot) \in
\cH^1(\T)$, i.e. the smoothness is less than in the case of the
Schr{\"o}dinger operator, where $\r_n(\cdot) \in \cH^2(\T)$. So we
need a specific version of the implicit function theorem (Theorem
\ref{p7t2}).

3) In the case of the Schr{\"o}dinger operator, each state runs
strictly monotonically around the gap $ \g_n^c$, changing sheets
when it hits the gaps boundary. In general, the states of the Dirac
operator on the half-line run non-monotonically and the direction of
their motion along $\g_n^c$ depends on the sign of $q_1 (t) +
\m_n(t)$, which is seen from (\ref{p1e5}) (see examples of monotone
and non-monotone motion on Fig.~2). Moreover, we show in
Theorem \ref{p1c2} that there exists a potential such
that a state runs non-monotonically in each open gap.

4) Theorem \ref{p0t1}, iv) gives that the number of the points $t_0
\in [0,1)$ such that $\m_n(t_0) = \a_n^{\pm}$ equals the number of
zeros of the first component of 2-periodic eigenfunction on the
interval $[0,1)$. This statement also holds in the case of the
Schr{\"o}dinger operator. Moreover, states of the Schr{\"o}dinger
operator run strictly monotonically around the gaps. Then the number
of their revolutions equals the number of points $t_0 \in [0,1)$
such that $\r_n(t_0) = a_n^{\pm}$. Since states of the Dirac
operator do not run strictly monotonically around the gaps, we
cannot associate the number of revolutions with the number of points
$t_0 \in [0,1)$ such that $\m_n(t_0) = \a_n^{\pm}$. In addition, it
follows from the oscillation theory for Schr{\"o}dinger operator
that $n$-th 2-periodic eigenfunction has exactly $n$ simple zeros in
$[0,1)$ (see e.g. \cite{LevSar}). But in general, this does not hold
in the case of the Dirac operator (see e.g. \cite{Tes}).
    \end{remark}
    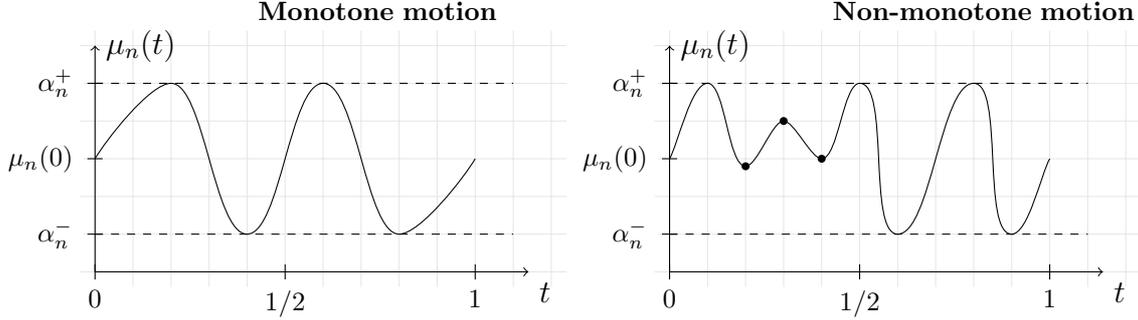
\begin{figure}[h]
            \begin{center}
                \begin{tikzpicture}
                    \usepgflibrary{decorations.markings}
                    \usepgfmodule{plot}
                    \usetikzlibrary{intersections}
                    \usetikzlibrary{decorations}
                    \usetikzlibrary{decorations.pathmorphing}
                    \usetikzlibrary{decorations.markings}
                    \usetikzlibrary{decorations.pathreplacing}
                    \usetikzlibrary{decorations.shapes}
                    \usetikzlibrary{topaths}
                    \usetikzlibrary{calc}
                    \usetikzlibrary{math}
                    \begin{scope}
                        \draw[ultra thin,gray!20,step=0.5cm] (-0.2,-0.2) grid (6.2,3.2);
                        \draw[->] (0,0) -- (5.7,0) node[anchor=north west] {$t$} coordinate(x axis);
                        \draw[->] (0,0) -- (0,3) node[anchor=west] {$\mu_n(t)$} coordinate(y axis);
                        \draw (0,.1) -- (0,-0.1) node[below,fill=white,font=\footnotesize] {$0$};
                        \draw (2.5,.1) -- (2.5,-0.1) node[below,fill=white,font=\footnotesize] {$1/2$};
                        \draw (5,.1) -- (5,-0.1) node[below,fill=white,font=\footnotesize] {$1$};
                        \draw (-0.1,0.5) -- (0.1,0.5) node[left=7pt,fill=white,font=\footnotesize] {$\alpha_n^-$};
                        \draw (-0.1,2.5) -- (0.1,2.5) node[left=7pt,fill=white,font=\footnotesize] {$\alpha_n^+$};
                        \draw (-0.1,1.5) -- (0.1,1.5) node[left=7pt,fill=white,font=\footnotesize] {$\mu_n(0)$};
                        \draw[dashed] (0,0.5) -- (5.5,0.5);
                        \draw[dashed] (0,2.5) -- (5.5,2.5);
                        \draw (0,1.5) to[in=180, out=60, looseness=0.5] (1,2.5) to[in=180, out=0, looseness=0.5] (2,0.5) to[in=180, out=0, looseness=0.5] (3,2.5) to[in=180, out=0, looseness=0.5] (4,0.5) to[in=-120, out=0, looseness=0.5] (5,1.5);
                        \draw (2,3.2) node[anchor = south west, text width=6cm, fill = white,font=\footnotesize]{\bf Monotone motion};
                    \end{scope}
                    \begin{scope}[xshift = 215]
                        \draw[ultra thin,gray!20,step=0.5cm] (-0.2,-0.2) grid (6.2,3.2);
                        \draw[->] (0,0) -- (5.7,0) node[anchor=north west] {$t$} coordinate(x axis);
                        \draw[->] (0,0) -- (0,3) node[anchor=west] {$\mu_n(t)$} coordinate(y axis);
                        \draw (0,.1) -- (0,-0.1) node[below,fill=white,font=\footnotesize] {$0$};
                        \draw (2.5,.1) -- (2.5,-0.1) node[below,fill=white,font=\footnotesize] {$1/2$};
                        \draw (5,.1) -- (5,-0.1) node[below,fill=white,font=\footnotesize] {$1$};
                        \draw (-0.1,0.5) -- (0.1,0.5) node[left=7pt,fill=white,font=\footnotesize] {$\alpha_n^-$};
                        \draw (-0.1,2.5) -- (0.1,2.5) node[left=7pt,fill=white,font=\footnotesize] {$\alpha_n^+$};
                        \draw (-0.1,1.5) -- (0.1,1.5) node[left=7pt,fill=white,font=\footnotesize] {$\mu_n(0)$};
                        \draw[dashed] (0,0.5) -- (5.5,0.5);
                        \draw[dashed] (0,2.5) -- (5.5,2.5);
                        \draw (0,1.5) to[in=180, out=60, looseness=0.5] (0.5,2.5) to[in=180, out=0, looseness=0.5] (1,1.4) node[fill=black, circle,minimum size=3pt,inner sep=0pt]{} to[in=180, out=0, looseness=0.5] (1.5,2) node[fill=black, circle,minimum size=3pt,inner sep=0pt]{} to[in=180, out=0, looseness=0.5] (2,1.5) node[fill=black, circle,minimum size=3pt,inner sep=0pt]{} to[in=180, out=0, looseness=0.5] (2.5,2.5) to[in=180, out=0, looseness=0.5] (3,0.5) to[in=180, out=0, looseness=0.5] (4,2.5) to[in=180, out=0, looseness=0.5] (4.5,0.5) to[in=-120, out=0, looseness=0.5] (5,1.5);
                        \draw (2,3.2) node[anchor = south west, text width=6cm, fill = white, font=\footnotesize]{\bf Non-monotone motion};
                    \end{scope}
                \end{tikzpicture}
            \label{fig2}
            \caption{Examples of monotone and non-monotone motion of $\m_n(t)$ in some open gap $\g_n$. The dots denote points where $\m_n(t)$ is not monotone.}
            \end{center}
        \end{figure}

In order to illustrate asymptotics (\ref{p1e3}--6), we consider the
case of continuous potential $q_1$. The following corollary show
that $\m_n(\cdot)$ moves locally as a state of the Schr{\"o}dinger
operator, when $q_1$ is continuous.

    \begin{corollary} \label{p1c1}
Let $V \in \cP$, $q_1 \in C(-\ve,\ve)$ for some $\ve > 0$, and let
a gap $\g_n$ be open for some $n \in \Z$ and $C_n = q_1(0) + \m_n^0 \neq
0$. Then $\m_n \in C^1(-\ve,\ve)$ and

\no i) If $\m_n^0 \neq \a_n^{\pm}$, then $\m_n$ is strictly monotone on $(-\tilde{\ve}, \tilde{\ve})$ for some $\tilde{\ve} > 0$ and satisfies as $t \to 0$:
\[ \label{p1e23}
\m_n(t) = \m_n^0 + \Omega(\m_n^0, 0) C_n t \left( 1 + o(1) \right).
\]
ii)  If $\m_n^0 = \a_n^{\pm}$, then $\m_n$ changes sheets at $t = 0$
and satisfies as $t \to 0$:
\begin{align}
\m_n(t) = \m_n^0 + \Omega_n^{\pm}(0) (C_n t)^2 \left(1 + o(1)
\right), \label{p1e24}
\\
(-1)^n\vp_2(1,\m_n(t),t) = 1 + \omega_n^{\pm}(0) C_n t \left(1 +
o(1)\right). \label{p1e26}
\end{align}
\end{corollary}
Remarks after Theorem \ref{p0t1} show  significant differences in
behavior of the states in the cases of the Dirac operator and of the
Schr{\"o}dinger operator. However, in Corollary \ref{p1c1} we show
that $\m_n(t)$ locally behaves as a state of the Schr{\"o}dinger
operator if $q_1$ is continuous. Now we prove that $\m_n(t)$
globally moves as a state of the Schr{\"o}dinger operator under some
restriction on the potential $q_1$.

\begin{theorem} \label{p0t2}
Let $V \in \cP$ and let a gap $\g_n$ be open for some $n \in \Z$.
Suppose that
\[
\label{p1e9}
\sign(q_1(t) + \a_n^-) = \sign(q_1(t) + \a_n^+) = \const \neq 0\ \
\text{for almost all $t \in [0,1)$.}
\]
Then the state $\m_n(\cdot)$ runs strictly monotonically around the
gap $\g_n^c$, changing sheets when it hits $\a_n^+$ or $\a_n^-$ and
making $\left| n \right|$ complete revolutions when $t$ runs through
$[0,1]$.
    \end{theorem}

Theorem \ref{p0t2} shows that under restrictions  (\ref{p1e9}) on
the potential $q_1$ the state $\m_n(t)$ moves as in the case of the
Schr{\"o}dinger operator.

\begin{corollary} \label{p0t4}
Let $V \in \cP$ and let a gap $\g_n$ be open for some $n \in
\Z$. Let one of the following conditions hold true:
        \begin{enumerate}[i)]
        \item $q_1 = 0$;
        \item $q_1 \in L^{\iy}(\T)$, $\left\|q_1\right\|_{\iy} < \min
        ( \left| \a_n^- \right|,\left| \a_n^+ \right|)$, and $0 \notin \g_n$;
            \item $q_1(t) > -\a_n^-$ for almost all $t \in [0,1)$;
            \item $q_1(t) < -\a_n^+$ for almost all $t \in [0,1)$;
        \end{enumerate}
Then the state $\m_n(\cdot)$ runs strictly monotonically around the
gap $\g_n^c$, changing sheets when it hits $\a_n^+$ or $\a_n^-$ and
making $\left| n \right|$ complete revolutions when $t$ runs through
$[0,1]$.
    \end{corollary}
    \begin{remark}
In the proof  we check that condition (\ref{p1e9}) holds true. If $V
= q_2 J_2$, then the Dirac operator is a supersymmetric charge (see
\cite{Thal}). In case ii), it is possible to show that the behavior
of the states of the Dirac operator does not differ from the
behavior of the states of the Schr{\"o}dinger operator at high
energies.
\end{remark}

Above we show, that $\m_n$ runs strictly monotone locally or
globally under some restriction on the potential. But in general the
motion is not monotone.

\begin{theorem} \label{p1c2}
i) There exists a potential $V \in \cP$ such that $\m_n(t)$ is not
monotone on $\T$ in each open gap $\g_n^c$, $n \in \Z$.

ii) There exists a potential $V \in \cP$ such that $\m_n(t)$ is not
monotone on $\T$ in any finite number of open gaps $\g_n^c$ and it
is monotone in infinite number of other open gaps.
\end{theorem}
\begin{remark}
    In the proof of this theorem we construct the needed potentials.
\end{remark}

\subsection{Tauberian theorem}

Now we show how one can use these results to construct a potential
from an asymptotic of $\m_n(t)$. In order to formulate these
results, we introduce a {\it slowly varying function}.
    \begin{definition*}
We say that $\zeta$ is slowly varying, if $\zeta$ is a positive
measurable function defined on some neighbourhood $(0,\ve)$, $\ve >
0$, and satisfying $\zeta(\l x)/\zeta(x) \to 1$ as $x \downarrow 0$,
for each $\l > 0$.
    \end{definition*}
    \begin{remark}
Usually the slowly varying functions are considered in a
neighborhood of infinity. In order to define such functions one need
replace $x$ by $1/x$ in the definition.
    \end{remark}
    \begin{definition*}
We say that $\xi(x) = x^{\varrho}\zeta(x)$ defined on $(0,\ve)$,
$\ve > 0$, is a regularly varying function if $\zeta$ is slowly
varying function on $(0,\ve)$ and $\varrho \in \R$.
    \end{definition*}
    \begin{remark}
1) In the definition $x^{\varrho}$ is positive for any $x \in
(0,\ve)$. It is easy to see that $\xi$ is regularly varying if and
only if $\xi(\l x)/ \xi(x) \to \l^{\varrho}$ as $x \downarrow 0$,
for each $\l > 0$.

2) In order to give a slight idea of the slowly and regularly
varying functions, we remark that $\const > 0$, $(\ln(1/x))^n$, $n
\in \R$, $\ln \ln (1/x)$, and similar are slowly varying but $x^n$
is not slowly varying for each $n \in \R$, $n \neq 0$. Thus, the
slowly varying functions decrease or increase at origin slowly than
$x^n$ for any $n \in \R$, $n \neq 0$. Remark that regularly varying
functions form an accurate asymptotic scale. These functions also
have many remarkable properties and they are widely used in the
formulation of Tauberian theorems (see e.g \cite{RegVar, Taub,
Seneta}).
    \end{remark}
    \begin{theorem} \label{p0t6}
Let $V \in \cP$ and let $q_1$ be a monotone function on neighborhood
$U^+ = (0,\ve)$ (or $U^- = (-\ve,0)$) for some $\ve > 0$. Let gap
$\g_n$ be open for some $n \in \Z$. Then the following asymptotic
holds true for some constants $C_{\pm} \in \R$, $\varrho_{\pm} > 0$,
and a slowly varying function $\zeta_{\pm}$:
\[ \label{p1e14}
\m_n(\pm t) - \m_n^0 = C_{\pm} t^{\varrho_{\pm}} \zeta_{\pm}(t) (1 +
o(1))\ \ \text{as $t \downarrow 0$}
        \]
if and only if the following asymptotics hold true for the same
constants $C_{\pm}$, $\varrho_{\pm}$, and slowly varying function
$\zeta_{\pm}$:
\begin{enumerate}[i)]
\item If $\m_n^0 \neq \a_n^{+}$, $\m_n^0 \neq \a_n^{-}$, then we have
\[ \label{p1e12}
q_1(\pm t) + \m_n^0 = \pm \frac{ C_{\pm}
\varrho_{\pm}}{\Omega(\m_n^0,0)} t^{\varrho_{\pm}-1} \zeta_{\pm}(t)
(1 + o(1))\ \ \text{as $t \downarrow 0$}.
 \]
 \item If $\m_n^0 = \a_n^{\kappa}$ for some $\kappa = \pm$ associated with the gap boundary, then we have
 \[ \label{p1e13}
 q_1(\pm t) + \m_n^0 = \mp \kappa S_n( \pm t) \frac{\left|C_{\pm}\right|^{1/2} \left|
  \varrho_{\pm} \right| }{2 \left| \Omega_n^{\kappa}(0) \right|^{1/2}} t^{\varrho_{\pm}/2-1}
  \zeta_{\pm}(t)^{1/2} (1 + o(1))\ \ \text{as $t \downarrow 0$},
 \]
 where $S_n(t) = \sign((-1)^n\vp_2(1,\m_n(t), t) - 1)$ and $\zeta_{\pm}(t)^{1/2} > 0$ for $t \in (0,\ve)$.
 \end{enumerate}
    \end{theorem}
\begin{remark}
1) In the theorem any real power of $\t \in (0,\ve)$ supposed to be
positive. We formulate the theorem in neighborhoods of $0$. Shifting
$V$ by any $t_0 \in \R$, one can obtain Theorem \ref{p0t6} in
a neighborhoods of $t_0$.

2) Identities (\ref{p1e12}--13) can be interpreted as derivatives of
asymptotic identities (\ref{p1e3}--5), where $\m_n(t)-\m_n^0$ is
determined from (\ref{p1e14}). From this point of view, one can
differentiate these identities because the conditions of theorem are
satisfied. See discussion and results about the problem of
differentiating an asymptotic expansion in \cite{Gra1}. In order to
prove the theorem, we use the Monotone Density Theorem (see Theorem
1.7.2b in \cite{RegVar}).

3) As we noted above, it follows from Lemma \ref{p3l1} that the
position of $\m_n( \pm t)$ on the Riemann surface is related to the
value of $S_n(\pm t)$. Thus, we can determine sign of $q_1(\pm t) +
\m_n^0$ in (\ref{p1e13}) by the position of $\m_n(\pm t)$ on the
Riemann surface.
\end{remark}
Using this theorem, one can construct asymptotic of $q_1$ from
asymptotic of $\m_n$. Remark that asymptotics of $q_1$ contains
$\Omega(\m_n^0,0)$ and $\left| \Omega_n^{\pm}(0) \right|$, so we can
determine asymptotics accurate to a constant factor. If we want to
construct a potential in each point we need to verify conditions of
Theorem \ref{p0t6}. So we need to check that $q_1$ is monotone and
there exist asymptotic (\ref{p1e14}) in neighborhood of each point.
We discuss this points in the following remark.
\begin{remark}
1) Using the definition of a compact set, it is easy to see that
$q_1$ is monotone in some left and right neighborhood of each point
if and only if it is piecewise monotone on the finite partition,
i.e. it has the form
$$
q_1(x) = \sum_{i = 1}^N q_{1i}(x) \chi_{(x_{i-1},x_i)}(x),\ \ x \in
[0,1] \sm \{x_i\}_{i=1}^N,
$$
where $\chi_{(a,b)}(x)$ is a characteristic function of the interval
$(a,b)$, $N > 0$, $q_{1i}$ is a monotone function for each $1 \leq i
\leq N$, and $0 = x_0 < x_1 < \ldots < x_{N-1} < x_N = 1$.

2) If $q_1$ has this form, then for each $t_0 \in [0,1]$ there
exists $\lim_{t \downarrow 0} q_1(t_0 \pm t)$, which may be
infinite. Using this fact, it is easy to see that if $C q_1(t_0 \pm
t)$ is a regularly varying function of $t$ for some constant $C \neq
0$ and $\lim_{t \downarrow 0} q_1(t_0 \pm t) \neq -\m_n(t_0)$, then
asymptotics (\ref{p1e12}--13) hold true for some $C_{\pm}$,
$\varrho_{\pm}$, and $\zeta_{\pm}$, which yields (\ref{p1e14}). In
order to show that the condition $\lim_{t \downarrow 0} q_1(t_0 \pm
t) \neq -\m_n(t_0)$ is crucial we consider the following simple
example. Let $q_1 \in L^2(\T)$ have the asymptotic expansion
$q_1(t_0 + t) = C + e^{-1/t}$ as $t \downarrow 0$, where $t_0 \in
\R$, $C \in \R$, $C \neq 0$, and $\m_n(t_0) = -C$. Thus, $q_1(t_0 +
\cdot)$ is slowly varying but $q_1(t_0+\cdot)+\m_n(t_0)$ is not
regularly varying.
\end{remark}

\subsection{Inverse problem}
As the application of Theorem \ref{p0t6}, we consider the inverse
problem for the potential $q_1 \in L^2(\T)$ defined by
\[
\label{p1e18}
q_1(t) = \sum_{i=1}^M \frac{C_i}{\left| t - t_i
\right|^{\varrho_i}} + \sum_{i=M+1}^{M+N} \frac{C_i
\sign(t-t_i)}{\left| t - t_i \right|^{\varrho_i}} +
\sum_{i=M+N+1}^{M+N+K} C_i (\sign(t - t_i) + D_i) ,\ \ t \in [0,1),
\]
where $M,N,K \geq 0$, $t_i \in [0,1)$ are distinct points for $1
\leq i \leq N+M+K$, and
$$
        \begin{aligned}
            &C_i \in \R \sm \{0\},\, \varrho_i \in (0,1/2), &&\text{for $1 \leq i \leq N+M$};\\
            &C_i \in \R \sm \{0\},\, D_i \in \R, &&\text{for $N+M+1 \leq i \leq N+M+K$}.
        \end{aligned}
    $$
\begin{theorem} \label{p0c1}
Let $V \in \cP$ and let $q_1$ have form (\ref{p1e18}). Let a gap
$\g_n$ be open for some $n \in \Z$. Then for each $t_0 \in [0,1]$
the following asymptotic holds true:
\[ \label{p1e19}
            \m_n(t_0 \pm t) - \m_n(t_0) = C_{\pm} t^{\varrho_{\pm}} (1 + o(1))\ \ \text{as $t \to 0$},
\]
for some $C_{\pm} \in \R$, and $\varrho_{\pm} > 1/2$. Moreover, if
for each $t_0 \in [0,1]$ asymptotic (\ref{p1e19}) and position of
$\m_n(t_0 \pm t)$ on the Riemann surface are given, then we can
recover all points $t_i$ and degrees $\varrho_i$. Furthermore, if we
know all constants $C_i$, $D_i$, then we can recover $q_1$.
\end{theorem}

\begin{remark}
1) In the proof of Theorem \ref{p0c1} we show how to recover
$t_i$, and explicitly express $\varrho_i$ through $\varrho_{\pm}$.

2) We consider a simple example to show how one can recover $q_1$
using Theorem \ref{p0c1}. Let $V \in \cP$ and let $q_1$ have form
(\ref{p1e18}). Let gap $\g_n$ be open and the following asymptotic
hold:
\[ \label{p1e22}
            \m_n(t_0 \pm \t) - \m_n(t_0) = C_0 \t^{1 - \varrho_0}(1 + o(1))\ \ \text{as $\t \downarrow 0$},
\]
for some $t_0 \in [0,1)$, $\m_n(t_0) \neq \a_n^{\pm}$, and $C_0 \in
\R \sm \{0\}$, $\varrho_0 \in (0,1/2)$. Moreover, let for any $t \in
[0,1) \sm \{t_0\}$ the following asymptotic hold:
$$
            \m_n(t \pm \t) - \m_n(t) = C_t (\pm \t)^{\varrho_t}(1 + o(1))\ \ \text{as $\t \downarrow 0$},
$$
for some $C_t \in \R \sm \{0\}$, and $\varrho_t \in \N$. Then the
application of Theorem \ref{p0c1} yields
        \[ \label{p1e21}
            q_1(t) = \frac{C_0(1-\varrho_0)}{\Omega(\m_n(t_0),t_0)}\frac{\sign(t-t_0)}{\left|t-t_0\right|^{\varrho_0}},\ \ t \in [0,1).
        \]
We used formula (\ref{p1e12}) to recover the potential. A more
detailed description is given in the proof of Theorem \ref{p0c1}.

3) Remark that we can use the asymptotic of $\m_n(t)$ in any open
gap to recover the potential. Moreover, if $q_1$ has form
(\ref{p1e21}), then in each open gap $\m_n(t)$ has asymptotic
(\ref{p1e22}) in neighborhood of $t_0$, so that $\m_n(t)$ is a
non-monotone function in each open gap.
    \end{remark}

\section{Periodic Dirac operator} \label{p2}
    \subsection{Dirac equation}
    We introduce the matrix-valued fundamental solution $\p(x, \l)$ of the Dirac equation
    \[ \label{p2e1}
        J y'(x) + V(x) y(x) = \l y(x),\ \ (x,\, \l) \in \R \ts \C.
    \]
satisfying the initial condition $\p(0,\l) = I_2$, where $I_2$ is
the identity matrix $2 \times 2$. Any matrix solution of some
initial value problem for equation (\ref{p2e1}) is expressed in
terms of $\p$ by multiplying on the right by the initial data.
$\p(1,\l)$ is called the monodromy matrix and satisfied
    \[ \label{p2e2}
        \p(1+x,\l) = \p(x,\l) \p(1,\l),\ \  (x,\, \l) \in \R \ts \C.
    \]
We introduce the vector-valued fundamental solutions $\vt(x,\l)$ and
$\vp(x,\l)$ of equation (\ref{p2e1}) satisfying the initial
conditions
    \[ \label{p2e3}
        \vt(0,\l) = \left( \begin{array}{c} 1 \\ 0 \end{array} \right),\ \
        \vp(0,\l) = \left( \begin{array}{c} 0 \\ 1 \end{array} \right).
    \]
    Moreover, $\p$ is expressed in terms of $\vp$ and $\vt$ as follows $\p = \left( \vt, \vp \right)$. We define the Wronskian of two vector-valued functions $f$ and $g$ by
    $$
        W(f,g) = f_1 g_2 - f_2 g_1.
    $$
It is known (see e.g. \cite{LevSar}) that Wronskian of $\vt$ and
$\vp$ does not depend at $x$ and satisfies
    \[ \label{p2e4}
        W(\vt, \vp) = \det \p = \vt_1 \vp_2 - \vp_1\vt_2 = 1.
    \]
We recall the known results about about $\p$ in the following
theorem (see e.g. \cite{Kor01}).

\begin{theorem} \label{p2t1}
Let $V \in \cP$. Then for each $\l \in \C$ there exists a unique
solution $\p$ of equation (\ref{p2e1}). For each $x \in [0,1]$ the
function $\p(x, \cdot)$ is entire. For each $\l \in \C$ the function
$\p(\cdot, \l), \p'(\cdot, \l) \in L^2(\T)$ and it satisfies
\[
\label{p2e8}
\left|\p(x,\l)\right| \leq e^{\left\|V\right\|_{\cP} +
\left|\Im \l\right| x},\ \ (x,\l) \in [0,1] \ts \C.
\]
    \end{theorem}

We introduce the Dirichlet eigenvalues $\m_n$, $n \in \Z$, as
eigenvalues of the Dirac equation (\ref{p1e1}) with boundary
conditions $y_1(0) = y_1(1) = 0$ and the Neumann eigenvalues $\n_n$,
$n \in \Z$, as eigenvalues of the Dirac equation (\ref{p1e1}) with
boundary conditions $y_2(0) = y_2(1) = 0$, where $y = \ma y_1 & y_2
\am^{\top}$. It is well known (see e.g. \cite{LevSar}) that $\m_n,
\n_n \in [\a_n^{-},\a_n^+]$, $n \in \Z$. If the gap $\g_n$
degenerates, then we have $\m_n = \n_n$. It follows from
(\ref{p2e3}) that eigenvalues $\m_n$ are zeros of an entire function
$\vp_1(1,\cdot)$, i.e. $\vp_1(1, \m_n) = 0$, $n \in \Z$. Thus,
$\vp(x,\m_n)$ is the eigenfunction for the eigenvalue $\m_n$ for
each $n \in \Z$. Below we need the following identity (see e.g.
\cite{Kor01})
    \[ \label{p2e19}
        \left\|\vp(\cdot,\m_n)\right\|^2 = - \vp_2(1,\m_n) \left. \partial_{\l} \vp_1(1,\l) \right|_{\l = \m_n},
    \]
where $\partial_{\l} u = u'_{\l} = \frac{\partial u}{\partial \l}$,
and $\left\| \cdot \right\|$ was defined above. The identity
(\ref{p2e19}) implies that $\pa_{\l} \vp_1(1,\m_n) \neq 0$ and each
Dirichlet eigenvalue is simple.

\subsection{Periodic Dirac operator}
For the Dirac equation we introduce the Lyapunov function
    \[ \label{p2e9}
        \D(\l) = \frac{1}{2} \Tr \p(1,\l) = \frac{1}{2}(\vp_2(1,\l) + \vt_1(1,\l)), \ \ \l \in \C.
    \]
    Due to Theorem \ref{p2t1} the function $\D$ is entire and describes a periodic spectrum by:
    $$
        \D(\a_n^{\pm}) = (-1)^{n};\ \ \left|\D(\l)\right| \leq 1,\, \l \in \s_n; \ \ \left|\D(\l)\right| > 1,\, \l \in \g_n,\ \ n \in \Z.
    $$
    For the periodic Dirac operator we define the Weyl-Titchmarsh function
    \[ \label{p2e10}
        m_{\pm}(\l) = \frac{a(\l) \mp b(\l)}{\vp_1(1,\l)},\ \ \l \in \mathbb{C}_+,
    \]
    where the functions $a(\l)$ and $b(\l)$ are defined by
    \[ \label{p2e11}
        a(\l) = \frac{\vp_2(1,\l) - \vt_1(1,\l)}{2},\ \ b(\l) = \sqrt{\D(\l)^2 - 1},\ \ \l \in \C_+,
    \]
where $i\sqrt{\D(\l+i0)^2 - 1} < 0$ for $\l \in \s_0$. The function
$a(\l)$ is entire and $b(\l)$ admits an analytic continuation from
$\C_+$ onto the Riemann surface $\L$. Recall that $\L$ consist of
two sheets $\L_1 = \L_2 = \C \sm \s(H^+_0)$ glued together crosswise
along the bands. In a neighborhood of the points $\a_n^{\pm}$ we use
the map $\l = \a_n^{\pm} \mp z^2$, where $z$ is the local variable
on the Riemann surface $\L$ and $\l \in \L_1$ if $z > 0$; $\l \in
\L_2$ if $z < 0$.  In addition, the following identity holds in each
open circle gap $\g_n^c = \overline{\g}_n^{(1)} \cup
\overline{\g}_n^{(2)} \ss \L$:
    \[ \label{p2e23}
        b(\l) = (-1)^{n+j+1} \left| \D^2(\l) - 1 \right|^{1/2},\ \ \l \in \g_n^{(j)},\ \ j=1,2.
    \]
Due to (\ref{p2e10}), $m_{\pm}(\l)$ has a meromorphic continuation
from $\C_+$ onto the Riemann surface $\L$. We introduce the Bloch
solutions $\p^{\pm}$ of equation (\ref{p2e1}) given by
    \[ \label{p2e20}
        \p^{\pm}(x,\l) = \vt(x,\l) + m_{\pm}(\l) \vp(x,\l),\ \ (x, \l) \in \R \ts \C_+.
    \]
It follows from analytic properties of $m_{\pm}$ that these
solutions admit a meromorphic continuation from $\C_+$ onto the
Riemann surface $\L$. If $\l \in \L$ is not a pole of $m_{\pm}$,
then these solutions are quasiperiodic, i.e.
    \[ \label{p2e13}
        \p^{\pm}(x+1,\l) = e^{\pm i k(\l)} \p^{\pm}(x,\l),\ \ x \in \R,
    \]
where $k(\l)$ is quasimomentum defined by $\D(\l)=\cos
k(\l)$. One can introduce quasimomentum as conformal mapping (see
\cite{Mi, KK}). We introduce the effective masses $M_n^{\pm} = 1/\l''(\a_n^{\pm})$, where $\l(k)$ is the inverse function for $k(\l)$. In \cite{KK} it was shown that $M_n^{\pm} = -\D(\a_n^{\pm})\D'(\a_n^{\pm})$ and $\pm M_n^{\pm} > 0$. Remark that quasimomentum is real-valued on the
spectral bands and $k(\a_n^{\pm}) = \pi n$, $n \in \Z$. This implies
that if $\l = \a_n^{+}$ or $\a_n^{-}$ and if $\l$ is not a pole of
$m_{\pm}$, then $\p^{+}(\cdot,\l) = \p^{-}(\cdot,\l)$ is periodic or
antiperiodic solutions of the Dirac equation. For any $\l \in
(\a_n^+, \a_{n+1}^-)$, $n \in \Z$, the solutions
$\p^{\pm}(\cdot,\l)$ are linearly independent, uniformly bounded on
the line, and do not decrease at infinity. On the other hand, if $\l
\in \g_n^{(1)}$ is not a pole of $m_{\pm}$, then $\p^{\pm}(x,\l)$
decreases exponentially as $x \to \pm \iy$ and increase
exponentially as $x \to \mp \iy$. Hence they are the Weyl solutions
for equation (\ref{p2e1}) and $\p^{\pm}(x,\l) \in
L^2(\mathbb{R}_{\pm})$. If $\l = \m_n$, then $\vp(x,\l)$ and
$\vt(x,\l)$ are the Weyl solutions and there are identities similar
to (\ref{p2e13}) for $\vp$ and $\vt$ (see Lemma \ref{p3l1}).

    \begin{lemma} \label{p2l3}
    \begin{enumerate}[i)]
        \item In any open gap $\g_n^c$, $n \in \Z$, the following asymptotic holds true:
        \[ \label{p2e12}
            b(\a_n^{\pm} \mp z^2) = (-1)^n z \sqrt{2\left|M_n^{\pm}\right|} + O(z^3)\ \ \text{as } z \to 0.
        \]
        Moreover, if $\m_n = \a_n^{\pm}$, then we get
        \[ \label{p2e18}
            2 (-1)^n M_n^{\pm} = -\vt_2(1,\a_n^{\pm})\left\|\vp(\cdot,\a_n^{\pm})\right\|^2.
        \]
        \item For any $\l \in \L$ the following identities hold true:
        \[ \label{p2e22}
            \begin{aligned}
                a^2(\l) - b^2(\l) &= -\vp_1(1,\l) \vt_2(1,\l), \\
                m_+(\l) m_-(\l) &= - \frac{\vt_2(1,\l)}{\vp_1(1,\l)}.
            \end{aligned}
        \]
    \end{enumerate}
    \end{lemma}
    \begin{proof}
        Let $\l = \a_n^{\pm} \mp z^2$ and $z \to 0$. Then we have for $z \to 0$
        \[ \label{p3e4}
            \D(\a_n^{\pm} \mp z^2) = \D(\a_n^{\pm}) \mp \D'(\a_n^{\pm})z^2 + O(z^4).
        \]
        From the definition of $\D(\cdot)$ it follows that $\D(\a_n^{\pm}) = (-1)^n$ and $-(-1)^n\D'(\a_n^{\pm}) > 0$. Using these facts and substituting asymptotic (\ref{p3e4}) in (\ref{p2e23}), we get for $z \to 0$
        $$
            b(\a_n^{\pm} \mp z^2) = (-1)^n\sign z \sqrt{\mp 2 \D(\a_n^{\pm})\D'(\a_n^{\pm})z^2 + O(z^4)}
            = (-1)^n z \sqrt{2\left| M_n^{\pm}\right|} + O(z^3).
        $$

        Let $\m_n = \a_n^{\pm}$. Differentiating (\ref{p2e1}) by $\l$, it is easy to see that $\p_{\l}(\cdot,\l)$ is a solution of the Dirac equation with the same potential and with inhomogeneous term $\p(x,\l)$. The solution of such equation can be represented as follows (see proof of (\ref{p2e17}) in Lemma \ref{p6l1})
        \[ \label{p2e17}
            \p_{\l}(1,\l) = \p(1,\l)\int_0^1 \p^{-1}(t,\l) J \p(t,\l) dt.
        \]
        Using the definition of $\D(\l)$ and identity (\ref{p2e17}), we get
        $$
            \D'(\a_n^{\pm}) = \frac{1}{2}\Tr \p_{\l}(1,\a_n^{\pm}) = \frac{1}{2}\vt_2(1,\a_n^{\pm})\left\| \vp(\cdot,\a_n^{\pm})\right\|^2.
        $$
        Finally, it follows $M_n^{\pm} = -\D(\a_n^{\pm})\D'(\a_n^{\pm})$ and $\D(\a_n^{\pm}) = (-1)^n$ that
        $$
            \vt_2(1,\a_n^{\pm})\left\|\vp(\cdot,\a_n^{\pm})\right\|^2 = -2(-1)^n M_n^{\pm}.
        $$

        In order to prove identities (\ref{p2e22}), we use the definitions of $a$, $b$, and $\D$. Then we get
        $$
            \begin{aligned}
                a^2 - b^2 &= \frac{\left( \vp_2(1,\cdot) - \vt_1(1,\cdot) \right)^2}{4} - \frac{\left( \vp_2(1,\cdot) + \vt_1(1,\cdot) \right)^2}{4} + 1 \\
                &= 1 - \vp_2(1,\cdot) \vt_1(1,\cdot) = -\vp_1(1,\cdot) \vt_2(1,\cdot),
            \end{aligned}
        $$
        where (\ref{p2e4}) have been used. The last identity and (\ref{p2e10}) yield
        $$
            \begin{aligned}
                m_+ m_- = \frac{a - b}{\vp_1(1,\cdot)}\frac{a + b}{\vp_1(1,\cdot)} = \frac{a^2 - b^2}{\vp_1^2(1,\cdot)}
                = -\frac{\vp_1(1,\cdot) \vt_2(1,\cdot)}{\vp_1^2(1,\cdot)} = - \frac{\vt_2(1,\cdot)}{\vp_1(1,\cdot)}.
            \end{aligned}
        $$
    \end{proof}

Note that (\ref{p2e22}) implies that the Dirichlet or Neumann
eigenvalue $\l_0 = \a_n^{\pm}$ for some $n \in \Z$ if and only if
$a(\l_0) = 0$. The second identity in (\ref{p2e22}) allows us to
compare the analytical properties of $m_+$ and $m_-$ on the Riemann
surface, because the right-hand side of this identity is a
meromorphic function on $\C$.

    \subsection{Shifted Dirac equation}
    Recall that in (\ref{p1e1}) we have introduced the shifted Dirac equation
    \[ \label{p2e14}
        Jy'(x) + V(x+t)y(x) = \l y(x),\ \ (x, \l) \in \R \ts \C,
    \]
    where $t \in \R$ is the shift parameter. If $V \in \cP$, then for each $t \in \R$ we get equation (\ref{p2e1}) with the potential $V(\cdot+t) \in \cP$. Hence for this equation there are all objects introduced for (\ref{p2e1}). We add the dependence on $t$ to these objects if they are not constant functions of $t$.

    Note that if $y(x,\l)$ is a solution of equation (\ref{p2e1}), then $\widetilde{y}(x,\l,t) = y(x+t,\l)$ is a solution of equation (\ref{p2e14}). Thus the fundamental matrix-valued solution of equation (\ref{p2e14}) is simply expressed in terms of the solution of equation (\ref{p2e1})
    \[ \label{p2e15}
        \begin{aligned}
            \p(x,\l,t) & = \p(x+t,\l)\p^{-1}(t,\l),\\
            \p(1,\l,t) & = \p(t,\l)\p(1,\l)\p^{-1}(t,\l).\\
        \end{aligned}
    \]
    Using (\ref{p2e9}), (\ref{p2e15}), and the fact that the traces of similar matrices are equal, we get
    $$
        \begin{aligned}
            \D(\l,t) = \frac{1}{2} \Tr \left( \p(t,\l)\p(1,\l)\p^{-1}(t,\l) \right) = \frac{1}{2} \Tr \p(1,\l) = \D(\l,0).
        \end{aligned}
    $$
    It gives that $b(\l)$, $k(\l)$, $M_n^{\pm}$ also does not depend on $t$. So we do not write the argument $t$ of these functions.

    \subsection{Dirac operator on the half-line}
In Section \ref{p1} we introduced the Dirac operators
$H^{\pm}_t = J \frac{d}{dx}+V_t$ acting on $L^2(\R_{\pm}, \C^2)$
with the Dirichlet boundary condition $f_1(0) = 0$, $f = \ma
f_1 & f_2 \am^{\top}$, where  $(t,V) \in \R\ts \cP$. It follows from
Weyl-Titchmarsh theory for the Dirac operator that $H^{\pm}_t$ are
self-adjoint (see e.g. \cite{W87}). For each $\l \in \C_+$ the
resolvent $\left( H^{\pm}_t - \l \right)^{-1}$ is an integral
operator acting on $L^2(\R_{\pm}, \C^2)$ with kernel
    \[ \label{p3e7}
        R^{\pm}(x,y,\l,t) =
        \left\{ \begin{array}{l}
            \p^{\pm}(x,\l,t) \otimes \vp(y,\l,t),\ \ \left|x\right| > \left|y\right|, \\
            \vp(x,\l,t) \otimes \p^{\pm}(y,\l,t),\ \ \left|x\right| < \left|y\right|,
        \end{array} \right.
    \]
    where $x,y \in \R_{\pm}$, and $(\l,t) \in \C_+ \ts \R$. In this formula $\otimes$ means the tensor product of vectors. In particular in the case $\left|x\right| > \left|y\right|$, we have
    $$
        R^{\pm}(x,y,\l,t) =
        \ma
            \p^{\pm}_1(x,\l,t) \vp_{1}(y,\l,t) & \p^{\pm}_1(x,\l,t) \vp_{2}(y,\l,t) \\
            \p^{\pm}_2(x,\l,t) \vp_{1}(y,\l,t) & \p^{\pm}_2(x,\l,t) \vp_{2}(y,\l,t)
        \am.
    $$
\begin{lemma} \label{p3l2}
Let $\eta \in C_o^{\iy}(\R_{\pm})$, $\eta \neq 0$, and let $t \in \R$. Then $f_{\pm}(\l,t) = \left( \left(
H^{\pm}_t- \l \right)^{-1} \eta, \eta \right)$ admits a meromorphic continuation
from $\C_+ \ss \L_1$ onto the Riemann surface $\L$ as function of $\l$ and its poles coincide
with the poles of $m_{\pm}(\cdot,t)$, moreover, their multiplicities
coincide.
\end{lemma}
\begin{proof}
        It is easy to see that $f_{\pm}(\l,t)$ is an analytic function of $\l$ for any $(\l,t) \in \C_+ \ts \R$. It follows from (\ref{p3e7}) that
        $$
            f_{\pm}(\l,t) = \int_{\R_{\pm}} \int_{\R_{\pm}} \eta(x) R^{\pm}(x,y,\l,t) \eta(y) dx dy = \int_{\R_{\pm}} \int_{\R_{\pm}} R^{\pm}_{\eta}(x,y,\l,t) dx dy,
        $$
        where we introduce $R^{\pm}_{\eta}(x,y,\l,t) = \eta(x) R^{\pm}(x,y,\l,t) \eta(y)$. The solution $\vp(x,\cdot,t)$ is entire and the solution
        $\p^{\pm}(x,\cdot,t) = \vt(x,\cdot,t) +
        m_{\pm}(\cdot,t)\vp(x,\cdot,t)$ is analytic on $\L$ for any $x,t \in
        \R$. Thus $R^{\pm}_{\eta}(x,y,\cdot,t)$ admits a meromorphic continuation
        from $\C_+$ onto the Riemann surface $\L$ and its poles coincide
        with the poles of $m_{\pm}(\cdot,t)$ for any $x,y,t \in \R$, and $\eta \in C_o^{\iy}(\R_{\pm})$. For any $(\l,t) \in \L \ts \R$ such that $\l$ is not a pole of $m_{\pm}(\cdot,t)$, the function $R^{\pm}_{\eta}(x,y,\l,t)$ is a continuous function of $x,y$ with compact support. Thus $f_{\pm}(\l,t)$ is correctly defined for such $(\l,t)$, which yields that $f_{\pm}(\l,t)$ admits a meromorphic continuation from $\C_+$ onto $\L$ and its poles coincide with the poles of $m_{\pm}(\cdot,t)$, moreover, their multiplicities coincide.
    \end{proof}
    We describe the spectrum of the operator on the half-line. These results is well known.
    \begin{theorem} \label{p3t1}
        Let $V \in \cP$ and let $t \in \R$. Then the following statements hold true:
        \begin{enumerate}[i)]
            \item $\s (H^{\pm}_t) = \s_{ac} (H^{\pm}_t) \cup \s_{d} (H^{\pm}_t),\ \ \s_{ac} (H^{\pm}_t) \cap \s_{d} (H^{\pm}_t) = \es$,
            \item $\s_{ac} (H^{\pm}_t) = \s_{ac} (H^{+}_0) = \bigcup_{n \in \Z} \s_n$,
            \item $\s_{d} (H^{\pm}_t) \subset \{ \m_n(t) \}_{n \in \Z}$.
        \end{enumerate}
    \end{theorem}

    \begin{proof}
        By Lemma \ref{p3l2}, for each $\eta \in C_o^{\iy}(\R)$, $\eta \neq 0$ and $t \in \R$ the function $f_{\pm}(\cdot,t)$ is meromorphic on $\L$ and it has only isolated poles over $\{ \m_n(t) \}_{n \in \Z}$ which yields iii). Using (\ref{p2e11}), we get $\Im f_{\pm}(\l,t) \neq 0$, $\l \in \R$ if and only if $\left| \D(\l) \right| \geq 1$. Thus, Theorem XIII.20 from \cite{RS4} yields that singular continuous spectrum of $H^{\pm}_t$ is absent and $\s_{ac} (H^{\pm}_t) = \cup_{n \in \Z} \s_n$. Since any solution of the Dirac equation does not belong to $L^2(\R_{\pm},\C^2)$ if $\left| \D(\l) \right| \geq 1$, it follows that there are no eigenvalues embedded into the absolutely continuous spectrum which yields i), and ii).
    \end{proof}

    The following lemma provides a relation between resonances and eigenvalues of $H^{+}_t$ and $H^{-}_t$. In addition, we show that a Dirichlet eigenvalue in each open gap is a state of the operators $H^{\pm}_t$, and their positions on the Riemann surface $\L$ are determined by the value of $\vp_2(1,\l,t)$. Let $\l_* \in \L$ be the projection of $\l \in \L$ on the other sheet of $\L$.
    \begin{lemma} \label{p3l1}
        Let $V \in \cP$ and let $t \in \R$. Then the following statements hold true:
        \begin{enumerate}[i)]
            \item If $\l \in \L$ is a state of $H^{\pm}_t$, then $\l_*$ is a state of $H^{\mp}_t$.
            \item If $\l \in \L$ is a state of $H^{\pm}_t$, then $\vp_1(1,\l,t) = 0$ and $\l \in \g_n^c \neq \es$ for some $n \in \Z$.
            \item If $\l \in \g_n^c \neq \es$ for some $n \in \Z$ and $\vp_1(1,\l,t) = 0$, then $\l$ is a state of $H^{+}_t$ and
            \begin{itemize}
                \item if $\left|\vp_2(1,\l,t)\right| < 1$, then $\l \in \g_n^{(1)}$ is an eigenvalue of $H^{+}_t$;
                \item if $\left|\vp_2(1,\l,t)\right| > 1$, then $\l \in \g_n^{(2)}$ is a resonance of $H^{+}_t$;
                \item if $\left|\vp_2(1,\l,t)\right| = 1$, then $\l = \a_n^+$ or $\l = \a_n^-$ is a virtual state of $H^{+}_t$.
            \end{itemize}
        \end{enumerate}
    \end{lemma}
    \begin{proof}
        i) From Lemma \ref{p3l2}, it follows that the poles of the resolvent coincide with the poles of the Weyl-Titchmarsh function. Thus, in order to prove i) it is sufficient to show that $\l$ is a pole of $m_{\pm}(\cdot,t)$ if and only if $\l_*$ is a pole of $m_{\mp}(\cdot,t)$. Since $\vp_1(1,\cdot,t)$, and $a(\cdot,t)$ are entire functions, it follows that $\vp_1(1,\l_*,t) = \vp_1(1,\l,t)$, and $a(\l_*,t) = a(\l,t)$. It follows from the definition of $b$ that $b(\l_*) = -b(\l)$. Therefore, we have $m_{\pm}(\l_*,t) = m_{\mp}(\l,t)$, which proves the statement.

        ii) It follows from Lemma \ref{p3l2} that if $\l$ is a state of $H^{\pm}_t$, then it is a pole of $m_{\pm}(\cdot,t)$ which yields $\vp_1(1,\l,t) = 0$ and the projection of $\l$ on $\C$ coincides with $\m_n(t) \in \g_n$ for some $n \in \Z$. It gives $\l \in \g_n^c$ for some $n \in \Z$. If $\g_n^c$ is close, then $m_{\pm}(\cdot,t)$ has no pole at $\l = \a_n^+ = \a_n^-$ since $b(\cdot)$ is an analytic function of $z^2$ in a neighborhood of $\l$. Thus we get $\g_n^c \neq \es$.

        iii) Let $\vp_1(1,\l,t) = 0$. Then $\vp_2(1,\l,t)$ is an eigenvalue of the monodromy matrix $\p(1,\l,t)$ corresponds to the eigenvector $c = \ma 0 & 1 \am^{t}$ and $\p(x, \l,t) c = \vp(x,\l,t)$. Using formula (\ref{p2e2}), we have for any $m \in \Z$
        \[  \label{p3e3}
            \vp(x+m,\l,t) = \p(x,\l,t) \p(1,\l,t)^m c = \vp_2(1,\l,t)^m \vp(x,\l,t) .
        \]
        Let $\left|\vp_2(1,\l,t)\right| < 1$. Then it follows from (\ref{p3e3}) that $\vp(x+n, \l,t)$ is exponentially decreasing function of $n > 0$, and hence $\vp(\cdot,\l,t) \in L^2(\R_+,\C^2)$ is an eigenfunction of $H^{+}_t$ corresponding to the eigenvalue $\l \in \g_n^{(1)}$. Similarly, if $\left|\vp_2(1,\l,t)\right| > 1$, then $\vp(\cdot,\l,t) \in L^2(\R_-,\C^2)$ is an eigenfunction of $H^{-}_t$ corresponding to the eigenvalue $\l_* \in \g_n^{(1)}$. Then it follows from i) that $\l \in \g_n^{(1)}$ is a resonance of $H^{+}_t$.

        Let $\left|\vp_2(1,\l,t)\right| = 1$. Then $\vt_1(1,\l,t) = \vp_2(1,\l,t) = \pm 1$ which yields $\D(\l) = \pm 1$ and $\l = \a_n^{\kappa}$, where $\kappa = +$ or $-$. Let us show that $\l$ is a pole of $m_{\pm}(\cdot,t)$ on the $\L$, i.e. there is a nonzero term with negative degree in the Laurent series by local variable in a neighborhood of $\l$. Since $\g_n^c$ is open, then $\l$ is a branch point for $m_{\pm}(\cdot, t)$. Using (\ref{p2e12}) and the fact that $a(\cdot, t)$ and $\vp_1(1,\cdot,t)$ are entire functions, we get $m_{\pm}(\l + z^2, t) = \mp C/z + O(1)$ as $z \to 0$,  where $C = (-1)^n\sqrt{2\left| M_n^{\kappa} \right|}/\partial_{\l}\vp_1(1,\l,t) \neq 0$. Therefore $m_{\pm}(\cdot,t)$ has a pole at $\l$ and $\l$ is a virtual state of $H^{\pm}_t$.
    \end{proof}

\section{Proof of the main theorems} \label{p3}

    \begin{proof}[\bf{Proof of Theorem \ref{p0t1}}]
        From Lemma \ref{p3l1}, it follows that for any $t \in \R$ in each open gap $\g_n^c$, $n \in \Z$, there exists a unique state $\m_n(t)$ of the operator $H^{\pm}_t$ and $\vp_1(1,\m_n(t),t) = 0$.

        i) Let $t_0 \in [0,1)$ and let the gap $\g_n^c$ be open for some $n \in \Z$. The state $\m_n(t)$ is a solution of the equation $\vp_1(1,\m_n(t),t) = 0$ for $t$ in some neighborhood of $t_0$.  In proof we use the implicit function theorem \ref{p7t2}. From (\ref{p2e19}) it follows that ${\vp_1}_{\l}(1,\m_n(t),t) \neq 0$. By Lemma \ref{p6l1}, the conditions of Theorem \ref{p7t2} hold true for $\vp_1(1,\l,t)$ in some neighborhood of $(t_0,\m_n(t_0))$ and Theorem \ref{p7t2} yields that $\m_n(\cdot) \in \cH^1(t_0-\ve,t_0+\ve)$ for some $\ve > 0$. Using Lemma \ref{p7l2} and the fact $V$ is an 1-periodic function, we get $\m_n(\cdot) \in \cH^1(\T)$.

        Differentiating the identity $\vp_1(1,\m_n(t),t)=0$ by $t$ and using (\ref{p6e1}), we get
        \begin{align*}
            \dot{\m}_n(t) = \left. {- \frac{\dot{\vp}_1(1,\l,t)}{\partial_{\l} \vp_1(1,\l,t)} }\right|_{\l = \m_n(t)}
            = \frac{(q_1(t)+\m_n(t)) (\vp_2(1,\m_n(t),t) - \vt_1(1,\m_n(t),t))}{\partial_{\l} \vp_1(1,\m_n(t),t)}
        \end{align*}
        for almost all $t \in \T$. Since $\vp_1(1,\m_n(t),t) = 0$, it follows from (\ref{p2e4}) that $\vt_1(1,\m_n(t),t) = 1/\vp_2(1,\m_n(t),t)$. Using this formula and identity (\ref{p2e19}), we get the Dubrovin equation
        \[ \label{p3e2}
            \dot{\m}_n(t) = \frac{(q_1(t)+\m_n(t))(1 - \vp_2(1,\m_n(t),t)^2)}{\left\|\vp(\cdot,\m_n(t),t)\right\|^2 } = \left( q_1(t)+\m_n(t) \right) \Omega(\m_n(t), t).
        \]

        ii) Recall that $\m_n^0 = \m_n(0)$. We introduce the following notation $f(\t) = q_1(\t) + \m_n^0$. Integrating (\ref{p3e2}), we get
        \[ \label{p4e6}
            \m_n(t) - \m_n^0 = \Omega(\m_n^0, 0) \int_{0}^t f(\t) d\t + \eta_1(t),
        \]
        where
        $$
            \begin{aligned}
                \eta_1(t) = \int_{0}^t f(\t) (\Omega(\m_n(\t),\t) - \Omega(\m_n^0,0)) d\t
                + \int_{0}^t (\m_n(\t) - \m_n^0) \Omega(\m_n(\t),\t) d\t.
            \end{aligned}
        $$
        Let us estimate $\eta_1(t)$ as $t \to 0$. In the proof we consider all asymptotics when $t \to 0$. Recall that we define $Q_n(t)$, and $U_n(t)$ in (\ref{p1e11}). It is easy to see that $ Q_n(t) \leq  U_n(t) $ for any $(t,n) \in \R \ts \Z$. We also introduce the following function
        $$
            M_n(t) = \max_{\t \in [0,t]} \left| \m_n(\t) - \m_n^0 \right|,\ \ (t,n) \in \T \ts \Z.
        $$
        We consider the case of $t > 0$, the proof for $t < 0$ is similar. We introduce $\Omega_o(\l,t) = 2a(\l,t)/\partial_{\l}\vp_1(1,\l,t)$. It follows from point i) of this theorem that $\Omega_o(\m_n(t),t) = \Omega(\m_n(t),t)$ for any $(n,t) \in \Z \ts \R$. Using Lemma \ref{p6l1}, we get that there exist $h \in L^2(\T)$ such that $\left| \partial_t \partial_{\l} \Omega_o(\l,t)\right| \leq h(t)$ for any $(\l,t) \in \g_n^c \ts \T$. Integrating this estimate, we get that there exist $C > 0$ such that $\left| \partial_{\l}\Omega_o(\l,t) \right| \leq C$ for any $(\l,t) \in \g_n^c \ts \T$. Since $\m_n(t) \in \g_n^c$ for any $t \in \T$, the application of the mean value theorem yields
        \[ \label{p4e17}
            \max_{\t \in [0,t]} \left| \Omega_o(\m_n(\t),\t) - \Omega_o(\m_n^0,\t) \right| = O(M_n(t)).
        \]
        Differentiating $\Omega_o(\m_n^0,t)$ by $t$ and using Lemma \ref{p6l1}, it is easy to see that
        \[ \label{p4e13}
            \begin{aligned}
                \partial_{t} \Omega_o(\m_n^0,t) = F_1(t) f(t) + F_2(t) (q_1(t) - \m_n^0) + F_3(t) q_2(t) + F_4(t),
            \end{aligned}
        \]
        where we introduce
        $$
            F_1(t) = \frac{2\vt_2\partial_{\l}\vp_1 + 4a\partial_{\l}a}{(\partial_{\l} \vp_1)^2},\quad F_2(t) = -\frac{2\vp_1}{\partial_{\l} \vp_1},\quad F_3(t) = -\frac{4a}{\partial_{\l} \vp_1},\quad F_4(t) = \frac{4a^2}{(\partial_{\l} \vp_1)^2},
        $$
        and $\vp_1 = \vp_1(1,\m_n^0,t)$, $\vt_2 = \vt_2(1,\m_n^0,t)$, $a = a(\m_n^0,t)$. Thus $F_i(\m_n^0,t)$, $i = 1,2,3,4$, are absolutely continuous functions of $t$. Differentiating them by $t$ and using Lemma \ref{p6l1}, we get
        \[ \label{p4e14}
            \partial_{t} F_i(t) = O\left( \left| f(t) \right| + \left| q_1(t) - \m_n^0 \right| + \left| q_2(t) \right| + 1 \right),\ \ i = 1,2,3,4.
        \]
        Integrating (\ref{p4e14}) and using the H{\"o}lder inequality, we obtain
        \[ \label{p4e12}
            F_i(t) - F_i(0) = \int_{0}^t \partial_{\t} F_i(\t) d\t = O\left( U_n(t) \right).
        \]
        Substituting (\ref{p4e12}) in (\ref{p4e13}) and using the H{\"o}lder inequality, we have
        \[ \label{p4e16}
            \max_{\t \in [0,t]} \left| \Omega_o(\m_n^0,\t) - \Omega_o(\m_n^0,0) \right| = \left| \int_{0}^t \left| \partial_{\t} \Omega(\m_n^0,\t) \right| d\t \right| = O\left( U_n(t) + U_n(t)^2\right) = O\left( U_n(t)\right).
        \]
        Thus, it follows from (\ref{p4e17}), (\ref{p4e16}), and $\Omega(\m_n(\t),\t) = \Omega_o(\m_n(\t),\t)$ that
        \[ \label{p4e19}
            \max_{\t \in [0,t]} \left| \Omega(\m_n(\t),\t) - \Omega(\m_n^0,0) \right| = O\left( M_n(t) + U_n(t) \right).
        \]
        Now estimating the maximum of (\ref{p4e6}), we get
        $$
            M_n(t) = O\left( \left(Q_n(t) + \left| t \right| M_n(t)\right)\left(1 + M_n(t) + U_n(t)\right) \right),
        $$
        which yields
        \[ \label{p4e18}
            M_n(t) = O \left( Q_n(t) \right).
        \]
        Recall that $ Q_n(t) \leq  U_n(t) $. Using this inequality and substituting (\ref{p4e19}), and (\ref{p4e18}) in the definition of $\eta_1(t)$, we have $\eta_1(t) = O \left( Q_n(t) U_n(t) \right)$.

        iii) Let $\m_n^0 = \a_n^{\pm}$. Then we have $\vp_2(1,\m_n^0,0) = (-1)^n$, $a(\m_n^0,0) = 0$ and it follows that
        \[ \label{p4e11}
            \begin{aligned}
                \Omega(\m_n^0,0) = F_2(0) = F_3(0) = F_4(0) = 0,\quad F_1(0) = 2 \Omega_n^{\pm}(0).
            \end{aligned}
        \]
        Substituting (\ref{p4e11}), and (\ref{p4e12}) in (\ref{p4e13}), we get
        \[ \label{p4e22}
            \max_{\t \in [0,t]} \left| \Omega_o(\m_n^0,\t) - 2 \Omega_n^{\pm}(0) \int_{0}^\t f(s) ds \right| = O\left(U_n(t)^2\right).
        \]
        Now using (\ref{p4e17}) and (\ref{p4e22}) as well as for (\ref{p4e19}), we get
        \[ \label{p4e36}
            \begin{aligned}
                &\max_{\t \in [0,t]} \left| \Omega(\m_n(\t),\t) - 2 \Omega_n^{\pm}(0) \int_{0}^\t f(s) ds \right| = O(U_n(t)^2 + M_n(t)),\\
                &\max_{\t \in [0,t]} \left| \Omega(\m_n(\t),\t) \right| = O(Q_n(t) + U_n(t)^2 + M_n(t)).
            \end{aligned}
        \]
        Integrating equation (\ref{p3e2}), we get
        \[ \label{p4e21}
            \m_n(t) - \m_n^0 = 2 \Omega_n^{\pm}(0) \int_{0}^t f(\t) \int_{0}^\t f(s) ds d\t + \eta_2(t) = \Omega_n^{\pm}(0) \left( \int_{0}^t f(\t) d\t \right)^2 + \eta_2(t),
        \]
        where
        $$
            \begin{aligned}
                \eta_2(t) &= \int_{0}^t f(\t)\left( \Omega(\m_n(\t),\t) -2 \Omega_n^{\pm}(0) \int_{0}^{\t} f(s) ds \right) d\t
                + \int_{0}^t (\m_n(\t) - \m_n^0) \Omega(\m_n(\t),\t) d\t.
            \end{aligned}
        $$
        Using (\ref{p4e36}), we estimate the maximum of (\ref{p4e21}). So we get
        $$
            \begin{aligned}
                M_n(t) &= O\left(Q_n(t)^2\right) + O\left(Q_n(t)(U_n(t)^2 + M_n(t))\right)\\
                &+ O\left( \left| t \right| M_n(t) (Q_n(t) + U_n(t)^2  + M_n(t)) \right).
            \end{aligned}
        $$
        Since $Q_n(t) \leq U_n(t)$, we see that
        \[ \label{p4e23}
            M_n(t) = O\left( Q_n(t)U_n(t)\right).
        \]

        We determine asymptotic of $\vp_2(1,\m_n(t),t)$ as $t \to 0$. By Lemma \ref{p6l1}, $\vp_2(1,\m_n(t),t)$ satisfies the conditions of Lemma \ref{p7l3}, so that it is an absolutely continuous function of $t$. Thus, we have
        \[ \label{p4e31}
            \vp_2(1,\m_n(t),t) = (-1)^n + \int_{0}^t \partial_{\t}(\vp_2(1,\m_n(\t),\t)) d\t,
        \]
        where $\vp_2(1,\m_n^0,0) = (-1)^n$ since $\m_n^0 = \a_n^{\pm}$. Using (\ref{p6e1}), we get
        \[ \label{p4e32}
            \partial_{t}(\vp_2(1,\m_n(t),t)) = (q_1(t)+\m_n(t)) \Omega_1(\m_n(t),t),
        \]
        where
        $$
            \Omega_1(\m_n(t),t) = \vt_2(1,\m_n(t),t) + \Omega(\m_n(t),t)\partial_{\l}\vp_2(1,\m_n(t),t).
        $$
        Combining (\ref{p4e31}), and (\ref{p4e32}), we obtain
        $$
            \vp_2(1,\m_n(t),t) - (-1)^n = \vt_2(1,\m_n^0,0) \int_{0}^t f(\t) d\t + \eta_3(t),
        $$
        where
        $$
            \begin{aligned}
                \eta_3(t) &= \int_{0}^t f(\t) (\vt_2(1,\m_n(\t),\t) - \vt_2(1,\m_n^0,0))d\t \\
                &+ \int_{0}^t f(\t) \Omega(\m_n(\t),\t)\partial_{\l}\vp_2(1,\m_n(\t),\t) d\t
                + \int_{0}^t (\m_n(\t) - \m_n^0) \Omega_1(\m_n(\t),\t) d\t.
            \end{aligned}
        $$
        Using Lemma \ref{p6l1} as well as for (\ref{p4e19}), we get
        \[ \label{p4e37}
                \max_{\t \in [0,t]} \left| \vt_2(1,\m_n(\t),\t) - \vt_2(1,\m_n^0,0) \right| = O\left(M_n(t) + U_n(t)\right).
        \]
        Now using (\ref{p4e36}), (\ref{p4e23}), (\ref{p4e37}) and estimating continuous functions of $\t$: $\Omega_1(\m_n(\t),\t)$, $\partial_{\l}\vp_2(1,\m_n(\t),\t)$ by some constant, we get
        $$
\begin{aligned}
\eta_3(t) = O\left( Q_n(t)(M_n(t) + U_n(t) + Q_n(t)^2 + U_n(t)^2) +
\left| t \right| M_n(t) \right) = O\left( Q_n(t) U_n(t)\right).
\end{aligned}
        $$
        By (\ref{p2e18}), we get $\vt_2(1,\m_n^0,0) = (-1)^n\omega_n^{\pm}(0)$ if $\m_n^0 = \a_n^{\pm}$. Thus, we have proved (\ref{p1e6}). Now using this formula, we obtain
        \[ \label{p4e33}
            1 - \vp_2(1,\m_n(t),t)^2 = -2 (-1)^n \omega_n^{\pm}(0) \int_{0}^t (q_1(\t) + \m_n^0) d\t    + O(Q_n(t)U_n(t)).
        \]
        Using (\ref{p2e19}) it is easy to see that $\left\| \vp(\cdot, \m_n(t),t) \right\|^2$ is an absolutely continuous function of $t$. Then using (\ref{p2e19}), and Lemma \ref{p6l1}, we have
        \[ \label{p4e34}
            \left\| \vp(\cdot, \m_n(t),t) \right\|^2 = \left\| \vp(\cdot, \m_n^0,0) \right\|^2 + O(U_n(t)).
        \]
        Combining (\ref{p4e33}), and (\ref{p4e34}) and using $\Omega(\m_n(t),t) = \Omega_o(\m_n(t),t)$,  we get
        \[ \label{p4e35}
            \Omega(\m_n(t),t)  = 2 \Omega_n^{\pm}(0) \int_{0}^{t} (q_1(\t) + \m_n(\t)) d\t + O(Q_n(t)U_n(t)).
        \]
        Substituting (\ref{p4e23}), and (\ref{p4e35}) in the definition of $\eta_2(t)$, we get $\eta_2(t) = O\left( Q_n(t)^2 U_n(t)\right)$.

        iv) Let the gap $\g_n$ be open for some $n \in \Z$. Let $\l_0 = \a_n^{\pm}$ for any choice of the sign. And let $Y_1(t_0,\l_0,0) = 0$ for some $t_0 \in [0,1)$, where $Y(\cdot,\l_0,0)$ is a 2-periodic eigenfunction of equation (\ref{p2e14}) for $\l = \l_0$, and $t = 0$. Consider the function $y(x) = Y(x+t_0,\l_0, 0)$. It is obvious that $y_1(0) = Y_1(t_0,\l_0, 0) = 0$. Since $Y(\cdot,\l_0,0)$ is periodic or anti-periodic, it follows that $y_1(1) = Y_1(1+t_0,\l_0, 0) = \pm Y_1(t_0,\l,0) = 0$. In addition, $y$ is a solution of the Dirac equation
        $$
            Jy'(x) + V(x+t_0)y(x) = \l_0 y(x).
        $$
        From this it follows that $\m_n(t_0) = \l_0$. Conversely, let $\m_n(t_0) = \l_0$ for some $t_0 \in [0,1)$. Consider the function $y(x) = \vp(x-t_0,\l_0,t_0)$. Then $y_1(t_0) = 0$ and $y(x)$ is a 2-periodic solution of the equation $Jy'(x) + V(x)y(x) = \l_0 y(x)$. So it follows that $Y(x,\l_0,0) = y(x)$ and $Y_1(t_0,\l_0,0) = 0$.
    \end{proof}
    \begin{proof}[\bf{Proof of Corollary \ref{p1c1}}]
        Since $q_1 \in C(-\ve,\ve)$ for some $\ve > 0$, it follows from (\ref{p1e5}) that $\dot{\m}_n \in C(-\ve,\ve)$, which yields $\m_n \in C^1(-\ve,\ve)$. It is easy to see that $q_1(\t) + \m_n^0 = C_n + o(1)$ as $\t \to 0$, where $C_n = q_1(0) + \m_n^0$. Using this asymptotic, we get for $t \to 0$
        $$
            \begin{aligned}
                \int_0^t (q_1(\t) + \m_n^0) d\t = C_n t + o(t),\quad Q_n(t) = O(t),\quad U_n(t) = o\left( \left| t \right|^{1/2} \right).
            \end{aligned}
        $$
        Substituting this asymptotic in (\ref{p1e3}--6) we obtain (\ref{p1e23}--9). It follows from (\ref{p1e23}) that $\m_n$ is strictly monotone function in some neighborhood of $0$. By (\ref{p1e26}), $\vp_2(1,\m_n(t),t)$ is strictly monotone function of $t$, so that $\m_n(t)$ changes sheets at $t=0$ if $\m_n(0) = \a_n^{\pm}$.
    \end{proof}
    \begin{proof}[\bf{Proof of Theorem \ref{p0t2}}]
        Let the gap $\g_n^c$ be open for some $n \in \Z$ and suppose that
        $$
            \sign(q_1(t) + \a_n^-) = \sign(q_1(t) + \a_n^+) = C \neq 0\ \ \text{for almost all $t \in \T$,}
        $$
        where $C = \pm 1$. We introduce the following notation $\m_n^0 = \m_n(t_0)$, and $f(t) = q_1(t) + \m_n(t_0)$. Since $\m_n^0 \in \g_n^c$ for any $t_0 \in \T$, we get $\sign f(t) = C \neq 0$ for almost all $t \in \T$. It now follows that $\int_{t_0}^t f(\t) d\t$ is a strictly monotone function of $t$ and for each $t,t_0 \in \T$ we have
        \[ \label{p4e38}
            \left| \int_{t_0}^t f(\t) d\t \right| = \left| \int_{t_0}^t \left| f(\t) \right| d\t \right|.
        \]
        We denote the right-hand side of (\ref{p4e38}) by $Q_n(t_0,t)$. Shifting $V$ by some $t_0 \in \T$, it is easy to see that asymptotics (\ref{p1e3}--6) have the following form.
        \begin{enumerate}
            \item If $\m_n^0 \neq \a_n^{\pm}$ for some $t_0 \in \T$, then for $t \to t_0$ the following asymptotic holds:
            \[ \label{p4e41}
                \m_n(t) = \m_n^0 + \Omega(\m_n^0, t_0) \int_{t_0}^t f(\t) d\t + O\left(Q_n(t_0,t)U_n(t_0,t)\right).
            \]
            \item If $\m_n^0 = \a_n^{\pm}$ for some $t_0 \in \T$, then for $t \to t_0$ the following asymptotics hold:
                \begin{align}
                    \m_n(t) = \m_n^0 + \Omega_n^{\pm}(t_0) \left( \int_{t_0}^t f(\t) d\t \right)^2 + O\left(Q_n(t_0,t)^2U_n(t_0,t)\right), \label{p4e42} \\
                    (-1)^n\vp_2(1,\m_n(t),t) = 1 + \omega_n^{\pm}(t_0)  \int_{t_0}^t f(\t) d\t + O\left(Q_n(t_0,t)U_n(t_0,t)\right), \label{p4e43}
                \end{align}
        \end{enumerate}
        where for each $(n,t_0,t) \in \Z \times \R \times \R$ we introduce
        $$
            U_n(t,t_0) = \left| t -t_0 \right| + \left| t - t_0\right|^{1/2} \left| \int_{t_0}^t \left| V(\t) + \m_n^0I_2 \right|^2 d\t \right|^{1/2}.
        $$
        Substituting (\ref{p4e38}) in asymptotics (\ref{p4e41}), (\ref{p4e42}) and using the fact that $\int_{t_0}^t f(\t) d\t$ is a strictly monotone function of $t$, we get that $\m_n(t) - \m_n^0$ is a strictly monotone function of $t$ in some neighborhood of each $t_0$, where $\m_n^0 \neq \a_n^{\pm}$. Substituting (\ref{p4e38}) in asymptotics (\ref{p4e43}), we get that $\vp_2(1,\m_n(t),t) - (-1)^n$ is a strictly monotone function of $t$ in some neighborhood of each $t_0$, where $\m_n^0 = \a_n^{\pm}$. It now follows that $\m_n(t)$ moves monotonically around the gap $\g_n^c$ and changes a sheet at the boundary of the gap.

        From Theorem \ref{p0t1}, iv) and monotonicity of the motion of $\m_n(t)$, it follows that the number of revolutions of $\m_n(t)$ when $t$ runs through $[0,1)$ coincides with the number of zeros of $Y_1(\cdot,\a_n^{\pm},0)$ on the interval $[0,1)$. By Lemma \ref{p8l2}, if $\m_n(0) = \a_n^{\pm}$, then $\vp(\cdot, \m_n(0),0)$ is a 2-periodic eigenfunction and its first component has $\left|n\right|$ zeros on the interval $[0,1)$. By Lemma \ref{p8l3}, if $\m_n(0) \neq \a_n^{\pm}$, then the first component of the 2-periodic solution $\p^{\pm}(\cdot,\a_n^{\pm},0)$ has $\left|n\right|$ zeros on the interval $[0,1)$.
    \end{proof}
    \begin{proof}[\bf{Proof of Corollary \ref{p0t4}}]
        Let the gap $\g_n^c$ be open for some $n \in \Z$. In order to prove the theorem we check the assumptions of Theorem \ref{p0t2}.

        If $q_1 = 0$, then we have $\m_0 = \a_0^- = \a_0^+ = 0$. Since $\g_n^c$ is open, it follows that $n \neq 0$ and for all $t \in \T$ we get
        $$
            \sign(q_1(t) + \a_n^-) = \sign \a_n^- = \sign n = \sign \a_n^+ = \sign(q_1(t) + \a_n^+).
        $$
        Thus, the assumptions of Theorem \ref{p0t2} are satisfied.

        Let $q_1 \in L^{\iy}(\T)$, $\min ( \left| \a_n^- \right|,\left| \a_n^+ \right|) > \left\|q_1\right\|_{\iy}$, and $0 \notin \g_n^c$, which yields $\sign \a_n^- = \sign \a_n^+$. If $\sign \a_n^- = 1$, then we have
        $$
            q_1(t) + \a_n^{\pm} > \min ( \left| \a_n^- \right|,\left| \a_n^+ \right|) - \left\|q_1\right\|_{\iy} > 0
        $$
        for almost all $t \in \T$. If $\sign \a_n^- = -1$, then we have
        $$
            q_1(t) + \a_n^{\pm} < \left\|q_1\right\|_{\iy} - \min ( \left| \a_n^- \right|,\left| \a_n^+ \right|) < 0
        $$
        for almost all $t \in \T$. Thus, in both cases the assumptions of Theorem \ref{p0t2} are satisfied.

        Let $q_1(t) > -\a_n^-$ for almost all $t \in \T$. It follows from the inequality $\a_n^- \leq \a_n^+$, that $q_1(t) + \a_n^{\pm} >  0$ for almost all $t \in \T$. Thus, the assumptions of Theorem \ref{p0t2} are satisfied. The proof for $q_1(t) < -\a_n^+$ is similar.
    \end{proof}
\begin{proof}[\bf{Proof of Theorem \ref{p1c2}}]
i) We consider $V = q_1 J_1 \in \cP$, where $q_1(x) =
{\sign(x-t_o)\/ \left| x -t_o \right|^{1\/4}}$, $x \in [0,1),
t_o={1\/2}$. It is easy to see that for $t \to t_o$
\[ \label{p4e46}
\begin{aligned}
\int_{t_o}^t (q_1(\t) + \m_n(t_o)) d\t &= \frac{4}{3} \left| t-t_o \right|^{3/4}(1 + o(1)),
\\
Q_n(t_o,t) = \int_{t_o}^t \left| q_1(\t) + \m_n(t_o)\right| d\t &=
\frac{4}{3} \left| t-t_o \right|^{3/4}(1 + o(1)),
            \end{aligned}
        \]
where $Q_n$ have been introduced in the proof of Theorem \ref{p0t2}.
Let $\g_n^c$ be open for some $n \in \Z$. If $\m_n(t_o) \neq
\a_n^{\pm}$, then substituting (\ref{p4e46}) in (\ref{p4e41}), we
get for $t \to t_o$
 $$
 \m_n(t) - \m_n(t_o) = \frac{4}{3}\Omega(\m_n(t_o),t_o) \left| t-t_o \right|^{3/4}(1 + o(1)),
$$
which yields that $\m_n(t)$ is not monotone in some neighborhood of
$t_o$. If $\m_n(t_o) = \a_n^{\pm}$, then substituting (\ref{p4e46})
in (\ref{p4e43}), we get for $t \to t_o$
$$
(-1)^n\vp_2(1,\m_n(t),t) - 1 = \frac{4}{3}\omega_n^{\pm}(t_o) \left|
t-t_o \right|^{3/4}(1 + o(1)),
$$
which yields that $\m_n(t)$ does not change sheets at $t =t_o$ and
thus $\m_n(t)$ runs non-monotone at $\g_n^c$ for $t$ in some
neighborhood of $t_o$.

ii) Let $N \geq 0$, and let $V = q_1 J_1 \in \cP$, where $q_1(x) = p
\chi_{(t_o-\d,t_o)}(x) - p \chi_{(t_o,t_o+\d)}(x)$ for some $p,\d >
0$. We show that each gap in the spectrum of the Dirac operator with
such potential $V$ is open instead of the gap $(0,0)$. We introduce
$\psi_c(x,\l)$, $(x,\l) \in \R \ts \C$, as a solution of the Dirac
equation $J\psi_c(x,\l) + cJ_1\psi_c(x,\l) = \l \psi_c(x,\l)$
satisfying the initial condition $\psi_c(0,\l) = I_2$. Solving this
equation, we obtain
        \[ \label{p4e47}
            \psi_c(x,\l) = \ma \cos \omega x & -\sqrt{\frac{\l-c}{\l+c}} \sin \omega x \\
                            \sqrt{\frac{\l+c}{\l-c}} \sin \omega x & \cos \omega x \am, \ \ (x,\l) \in \R \ts \C,
        \]
        where $\omega = \sqrt{\l^2 -c^2}$. It is easy to see that the monodromy matrix for the Dirac operator with the potential $V$ has the following form
        \[ \label{p4e48}
            \psi(1,\l) = \psi_{0}(1/2-\d,\l) \psi_{-p}(\d,\l) \psi_{p}(\d,\l) \psi_{0}(1/2-\d,\l).
        \]
        Using (\ref{p4e47}), we obtain
        \[ \label{p4e49}
            \psi_{0}(1/2-\d,\l) = \ma c^0 & -s^0 \\ s^0 & c^0 \am,\ \ \psi_{-p}(\d,\l) = \ma c & -\frac{s}{m} \\ ms & c \am,\ \ \psi_{p}(\d,\l) = \ma c & -ms \\ \frac{s}{x} & c \am,
        \]
        where we introduce $c^0 = \cos \l(1/2-\d)$, $s^0 = \sin \l(1/2-\d)$, $c = \cos \sqrt{\l^2 -p^2}\d$, $s = \sin \sqrt{\l^2 -p^2}\d$, $m = \sqrt{\frac{\l+p}{\l-p}}$. Some gap $(\l,\l)$ is close if and only if $\psi(1,\l) = \pm I_2$. Substituting (\ref{p4e49}) in (\ref{p4e48}), we can rewrite equality $\vt_1(1,\l) = \vp_2(1,\l)$ as $\l-p = \pm(\l+p)$. Therefore if $p \neq 0$ then only the gap $(0,0)$ is close.

As in previous point of the theorem, using asymptotics
(\ref{p4e41}--24), one can prove that $\m_n(t)$ is not monotone at
$t = 1/2$ if $\max \{ \left| \a_n^+ \right|, \left| \a_n^- \right|
\} < \left| p \right|$. Since $V \in L^{\iy}(\T)$, there exist a
finite number of such gaps. Moreover, by Corollary \ref{p0t4},
$\m_n(t)$ is not monotone only in a finite number of gaps. Let us
show that there exist at least $N$ gaps $(\a^-_i, \a^+_i) \ss (-p,
p)$. Using Theorem 1.1 from \cite{Kor05b} and the definition of
$q_1$, we get for any $i_0 \in \Z$
\[ \label{p4e50}
            \sum_{i = i_0}^{i_0 + N} \left| \g_i \right| \leq \sqrt{N} \left( \sum_{i = i_0}^{i_0 + N} \left| \g_i \right|^2 \right)^{1/2} \leq \sqrt{N} \left( \sum_{i \in \Z} \left| \g_i \right|^2 \right)^{1/2} \leq 2\sqrt{2N} \left\| q_1 \right\| \leq 4p\sqrt{N\d},
        \]
where $\left| \g_n \right|$ is length of the gap. It follows from
Theorem 3.2 \cite{KK} that $\left| \s_n \right| \leq \pi$ for any $n
\in \Z$. Using this estimate, (\ref{p4e50}), and setting $p$, and
$\d$ such that $4\sqrt{N\d} < 1$, and $p > 4p\sqrt{N\d} + \pi N$, we
get for any $i_0 \in \Z$
$$
\sum_{i = i_0}^{i_0 + N} \left( \left| \g_i \right| + \left| \s_i
\right| \right) \leq A < p, \qqq A=4p\sqrt{N\d} + \pi N
        $$
        which yields that there exist at least $2N$ open gaps at the interval
        $
            (-A, A) \ss (-p,p).
        $
    \end{proof}
\begin{proof}[\bf{Proof of Theorem \ref{p0t6}}]
        Let $q_1$ be a monotone function on $U^+ = (0,\ve)$ for some $\ve > 0$. In the case when $q_1$ is a monotone function on $U^- = (-\ve,0)$, the proof is similar. Recall that $\m_n^0 = \m_n(0)$. We introduce the following notation $f(t) = q_1(t) + \m_n^0$. Since $q_1$ is monotone on $U^+$, we get $f(t)$ is also a monotone function on $U^+$. It now follows that there exist $\tilde{\ve}$ such that $0 < \tilde{\ve} \leq \ve$ and $\sign f(t) = \const$ for almost all $t \in (0,\tilde{\ve})$. Using this fact, we get for each $t \in (0,\tilde{\ve})$
        \[ \label{p4e44}
            \left| \int_{0}^t f(\t) d\t \right| = \left| \int_{0}^t \left| f(\t) \right| d\t \right| = \left| Q_n(t) \right|.
        \]
        In the proof we consider all asymptotics when $t \downarrow 0$. Using (\ref{p4e44}), we can rewrite asymptotic (\ref{p1e3}) for $\m_n^0 \neq \a_n^{\pm}$ as follows
        \[ \label{p4e24}
            \m_n(t) - \m_n^0 = \Omega(\m_n^0, 0) \int_{0}^t f(\t) d\t \left(1 + o(1)\right).
        \]
        If $\m_n^0 = \a_n^+$ or $\m_n^0 = \a_n^-$, then using asymptotics (\ref{p1e4}), (\ref{p1e6}), we have
        \begin{align}
            \m_n(t) - \m_n^0 = \Omega_n^{\pm}(0) \left( \int_{0}^t f(\t) d\t \right)^2 \left(1 + o(1)\right), \label{p4e25} \\
            (-1)^n \vp_2(1,\m_n(t),t) - 1 = \omega_n^{\pm}(0)  \int_{0}^t f(\t) d\t \left( 1 + o(1) \right). \label{p4e26}
        \end{align}
        Let asymptotics (\ref{p1e12}), (\ref{p1e13}) hold true for some $C_+ \in \R$, $\varrho_+ > 0$, and slowly varying function $\zeta_+$. Then integrating (\ref{p1e12}), (\ref{p1e13}) and using Proposition 1.5.8 from \cite{RegVar}, we get
        $$
            \begin{aligned}
                \int_{0}^{t} f(\t) d\t &= \frac{C_{+}}{\Omega(\m_n^0,0)} t^{\varrho_{+}} \zeta_{+}(t)(1+o(1)),\\
                \int_{0}^{t} f(\t) d\t &= - S_n(t) \frac{\left|C_{+}\right|^{1/2}}{\left| \Omega_n^{\kappa}(0) \right|^{1/2}} t^{\varrho_{+}/2} \zeta_{+}(t)^{1/2}(1+o(1)),
            \end{aligned}
        $$
        where $\m_n^0 = \a_n^{\kappa}$ in the second formula and $\kappa = \pm$ associated with the gap boundary. Remark that in the second case we get $\sign C_{+} = -\kappa$ since $\a_n^- \leq \m_n(t) \leq \a_n^+$ for any $(t,n) \in \R \ts \Z$. Substituting this result in (\ref{p4e24}), (\ref{p4e25}), we obtain (\ref{p1e14}).

        Let now asymptotic (\ref{p1e14}) holds true for some $C_+ \in \R$, $\varrho_+ > 0$, and slowly varying function $\zeta_+$.  If $\m_n^0 \neq \a_n^{\pm}$, then using (\ref{p4e24}), we have
        \[ \label{p4e27}
            \int_{0}^t f(\t) d\t = \frac{\m_n(t) - \m_n^0}{\Omega(\m_n^0, 0)} \left(1 + o(1)\right).
        \]
        Substituting (\ref{p1e14}) in (\ref{p4e27}), and using Theorem 1.7.2b from \cite{RegVar}, we get (\ref{p1e12}).

        If $\m_n^0 = \a_n^{\kappa}$ for some $\kappa = \pm$, then using (\ref{p4e25}), we have
        \[ \label{p4e28}
            \int_{0}^t \left| f(\t) \right| d\t = \frac{\left| \m_n(t) - \m_n^0 \right|^{1/2}}{\left| \Omega_n^{\kappa}(0) \right|^{1/2}} \left(1 + o(1)\right).
        \]
        Substituting (\ref{p1e14}) in (\ref{p4e28}), and using Theorem 1.7.2b from \cite{RegVar}, we get
        \[ \label{p4e29}
            \left| f(t) \right| = \frac{\left|C_{+}\right|^{1/2} \left| \varrho_{+} \right|}{2\left| \Omega_n^{\kappa}(0) \right|^{1/2}} t^{\varrho_{+}/2-1} \zeta_{+}(t)^{1/2} (1 + o(1)).
        \]
        We find $\sign f(t)$ from (\ref{p4e26}). Thus, we have for sufficiently small $t > 0$
        \[ \label{p4e30}
            \sign f(t) = S_n(t) \sign \omega_n^{\kappa}(0),
        \]
        where $S_n(t) = \sign((-1)^n \vp_2(1,\m_n(t), t) - 1)$, and $\sign \omega_n^{\kappa}(0) = - \sign M_n^{\kappa} = -\kappa$. Combining (\ref{p4e29}), and (\ref{p4e30}), we obtain (\ref{p1e13}).
    \end{proof}
    \begin{proof}[\bf{Proof of Theorem \ref{p0c1}}]

        Since $q_1$ has the form (\ref{p1e18}), it follows that $q_1$ is a piecewise rational function on some finite partition. The derivative of rational function is a rational function and has a finite number of zeros. It follows that there exist a finite number of point where $q_1$ change the monotonicity, which yields that $q_1$ is a monotone function on $U^+ = (t_0, t_0+\ve)$ and $U^- = (t_0- \ve, t_0)$ for some $\ve > 0$ and for each $t_0 \in [0,1]$.

        Using (\ref{p1e18}), we can distinguish four different types of the asymptotic behavior of $q_1(t_0 \pm \t) + \m_n(t_0)$ for $t_0 \in [0,1]$, $\t \downarrow 0$. Then the application of Theorem \ref{p0t6} yields that there exist four different types of the asymptotic behavior of $\m_n(t_0 \pm \t) - \m_n(t_0)$ and for each $t_0 \in [0,1]$
        $$
            \m_n(t_0 \pm \t) - \m_n(t_0) = C_{\pm} \t^{\varrho_{\pm}} (1 + o(1))\ \ \text{as $\t \downarrow 0$}.
        $$
        Describing these types, we show how to recover $t_i$ and $\varrho_i$. Let $\m_n(t_0) \neq \a_n^{\pm}$ for some $t_0 \in [0,1]$. Then we have the following types:
        \begin{enumerate}
            \item If $\varrho_{\pm} < 1$, and $\sign C_- \neq \sign C_+$, then we get $t_0 = t_i$ for some $1 \leq i \leq M$.
            \item If $\varrho_{\pm} < 1$, and $\sign C_- = \sign C_+$, then we get $t_0 = t_i$ for some $M+1 \leq i \leq M+N$.
            \item If $\varrho_{\pm} \in \N$, and $\varrho_+ \neq \varrho_-$ or $C_+ \neq (-1)^{\varrho_{\pm}} C_-$, then we get $t_0 = t_i$ for some\\ $M+N+1 \leq i \leq M+N+K$.
            \item If $\varrho_{\pm} \in \N$, $\varrho_+ = \varrho_-$, and $C_+ = (-1)^{\varrho_{\pm}} C_-$, then we get $t_0 \neq t_i$ for each\\ $1 \leq i \leq M+N+K$.
        \end{enumerate}
        Here $\N$ is the set of positive integer numbers. Moreover, we can recover $\varrho_i$ by the formula $\varrho_i = 1 - \varrho_{\pm}$. Let $\m_n(t_0) = \a_n^{+}$ or $\m_n(t_0) = \a_n^{-}$ for some $t_0 \in [0,1]$. Then we have the following types:
        \begin{enumerate}
            \item If $\varrho_{\pm} < 2$, and $S_n(t_0 + \ve) \neq S_n(t_0 - \ve)$ for all sufficiently small $\ve > 0$, then we get $t_0 = t_i$ for some $1 \leq i \leq M$.
            \item If $\varrho_{\pm} < 2$, and $S_n(t_0 + \ve) = S_n(t_0 - \ve)$ for all sufficiently small $\ve > 0$, then we get $t_0 = t_i$ for some $M+1 \leq i \leq M+N$.
            \item If $\varrho_{\pm} \in 2\N$, and $\varrho_+ \neq \varrho_-$ or $S_n(t_0 + \ve) \neq (-1)^{\varrho_{\pm}/2} S_n(t_0 - \ve)$ for all sufficiently small $\ve > 0$, then we get $t_0 = t_i$ for some $M+N+1 \leq i \leq M+N+K$.
            \item If $\varrho_{\pm} \in 2\N$, $\varrho_+ = \varrho_-$, and $S_n(t_0 + \ve) = (-1)^{\varrho_{\pm}/2} S_n(t_0 - \ve)$ for all sufficiently small $\ve > 0$, then we get $t_0 \neq t_i$ for each $1 \leq i \leq M+N+K$.
        \end{enumerate}
        Remark that $S_n(t_0 + \ve) \neq S_n(t_0 - \ve)$ for all sufficiently small $\ve > 0$ if and only if $\m_n(t)$ changes the sheet of the Riemann surface when $t$ runs through $(t_0 - \ve, t_0 + \ve)$. Moreover, we can recover $\varrho_i$ by the formula $\varrho_i = 1 - \varrho_{\pm}/2$. If we also know constants $C_i, D_i$, we can recover the potential using (\ref{p1e18}).
    \end{proof}

\section{Appendix} \label{p6}

    \subsection{Implicit function theorem}
In order to prove the main theorems we need a specific version of
the implicit function theorem. Recall that we have denoted the
Sobolev spaces on $\T$ by $\cH^i(\T)$, $i \geq 0$. We also denote
the Sobolev spaces on $I \ss \R$ by $\cH^i(I)$, $i \geq 0$. We
recall the standard form of the implicit function theorem (see
\cite{IFT}, \cite{Kum80}).

\begin{theorem}[Standard IFT] \label{p7t1}
        Let $F \in C(I \times I)$, where $I = [-1,1]$, and let
        \begin{enumerate}[i)]
            \item $F(0,0) = 0$,
            \item $F(x,\cdot)$ be strictly monotone on $I$ for each fixed $x \in I$.
        \end{enumerate}
Then there exist an open neighborhood $U_1 \times U_2 \ss I \times
I$ of $(0, 0)$ and a unique $f \in C(U_1)$ such that $f(U_1) \ss
U_2$ and for any point $(x,y) \in U_1 \times U_2$
        $$
            F(x,y) = 0 \text{ if and only if $y = f(x)$}.
        $$
    \end{theorem}
\begin{remark}
If $F \in C^1(I\ts I)$, then the condition of strict monotonicity
can be replaced by the condition $F_y(0,0) \neq 0$. Moreover, then
we get $f \in C^1(U_1)$.
\end{remark}
We need the implicit function theorem with weaker conditions.
\begin{definition*}
Let $I_1,I_2 \ss \R$ be open bounded intervals. We denote by
$\mH(I_1,I_2)$ the set of all functions $F:I_1 \times I_2 \to \R$
satisfying the following conditions
        \begin{enumerate}[i)]
            \item $F(x,\cdot) \in C^1(I_2)$ for each fixed $x \in I_1$,
            \item $F(\cdot,y) \in \cH^1(I_1)$ for each fixed $y \in I_2$,
            \item $\left|F'_{x} (x,y)\right| \leq g(x)$ for each $(x,y) \in I_1 \ts I_2$, for some $g \in L^2(I_1)$.
        \end{enumerate}
    \end{definition*}
    We give the following simple lemmas without proofs.
    \begin{lemma} \label{p7l1}
        We have $\mH(I,I) \ss C(I \ts I)$, where $I = [-1,1]$.
    \end{lemma}
    \begin{lemma} \label{p7l3}
        Let $F,F'_y \in \mH(I,I)$, where $I = [-1,1]$, and let $f \in \cH^1(I)$ such that $f(I) \ss I$. Then $F(\cdot,f(\cdot)) \in \cH^1(I)$.
    \end{lemma}
    \begin{lemma}\label{p7l2}
        Let $f: I \to \R$ be measurable, where $I = [-1,1]$. Let for each $x \in I$ there exists an open interval $I_x \ss I$ (in the subspace topology on $I$) such that $x \in I_x$, and $\left. f \right|_{I_x} \in \cH^1(I_x)$. Then we have $f \in \cH^1(I)$.
    \end{lemma}

    Let us prove now the specific form of the implicit function theorem.
    \begin{theorem}[Specific IFT] \label{p7t2}
        Let $F, F'_y \in \mH(I,I)$, where $I = [-1,1]$, and let
        \begin{enumerate}[i)]
            \item $F(0,0) = 0$,
            \item $F'_{y}(0,0) \neq 0$.
        \end{enumerate}
        Then there exist an open neighborhood $U_1 \times U_2 \ss I \times I$ of $(0, 0)$ and a unique $f \in \cH^1(U_1)$ such that $f(U_1) \ss U_2$ and for any point $(x,y) \in U_1 \times U_2$
        $$
            F(x,y) = 0 \text{ if and only if $y = f(x)$}.
        $$
    \end{theorem}
    \begin{proof}
        First we prove the existence of a continuous function $f$ and then we prove that it is absolutely continuous. The application of Lemma \ref{p7l1} yields $F,F'_y \in C(I \ts I)$. Since $F_y(0,0) \neq 0$ and $F'_y \in C(I \ts I)$, it follows that there exists some neighborhood $Y_1 \ts Y_2 \ss I \ts I$ of $(0,0)$, where $F'_y$ does not change sign and is separated from zero. So that $F(x,\cdot)$ is strictly monotone on $Y_2$ for each fixed $x \in Y_1$. Thus the conditions of Theorem \ref{p7t1} are satisfied and hence there exists an open neighborhood $U_1 \ts U_2$ of $(0,0)$ and a unique $f \in C(U_1)$ such that $f(U_1) \ss U_2$ and for any $(x,y) \in U_1 \ts U_2$
        \[ \label{p7e1}
            F(x,y) = 0\ \ \text{if and only if}\ \ y = f(x).
        \]
        Let $x_1,x_2 \in U_1$. Then using $F \in \mH(I,I)$, and (\ref{p7e1}), we get
        \[ \label{p7e4}
            \int_{x_1}^{x_2} F'_{x}(x,f(x_2)) dx = -\int_{f(x_1)}^{f(x_2)} F'_{y}(x_1,y) dy.
        \]
        Above we show that $F'_y$ is separated from zero on $U_1 \ts U_2 \ss Y_1 \ts Y_2$. Thus, there exists a constant $C > 0$ such that $\left|F'_y(x,y)\right| \geq C$ for any $(x,y) \in U_1 \ts U_2$. Using this inequality, we can evaluate the absolute value of the integral in the right-hand side of (\ref{p7e4}) from below. Thus, we have
        \[ \label{p7e3}
            \int_{x_1}^{x_2} \left|F'_{x}(x,f(x_2))\right| dx \geq C \left|f(x_2) - f(x_1)\right|.
        \]
        Since $F \in \mH(I,I)$, we can estimate $F'_x$ in (\ref{p7e3}). Thus, we get for any $x_1,x_2 \in U_1 \ss I$
        $$
            \left|f(x_2) - f(x_1)\right| \leq \frac{1}{C} \int_{x_1}^{x_2} g(x) dx,
        $$
        where $g \in L^2(U_1) \ss L^1(U_1)$, which yields that $f$ is absolutely continuous on $U_1$. It remains to show that $f' \in L^2(U_1)$. Differentiating the identity $F(x,f(x))=0$ by $x$, we get
        $$
            f'(x) = -\left. \frac{ F'_{x}(x,y) }{F'_{y}(x,y)} \right|_{y = f(x)},\ \ x \in U_1.
        $$
        Since $\left|F'_{y}(x,y)\right| \geq C > 0$ and $\left|F'_x(x,y)\right| \leq g(x)$ for any $(x,y) \in U_1 \ts U_2$, it follows that $\left|f'(x)\right| \leq g(x)/C$, $x \in U_1$, which yields $f' \in L^2(U_1)$ since $g \in L^2(U_1)$.
    \end{proof}
    \begin{remark}
        Each implicit function theorem corresponds to an existence (and uniqueness) theorem for differential equations (see e.g. \cite{abIFT}). Theorem \ref{p7t2} in some sense associated with the Carath{\'e}odory existence theorem (see e.g. Theorem 5.1 in \cite{Hall}), which says that there exists an absolutely continuous solution of the equation $\dot{f}(t) = G(t,f(t))$ through $(0,0)$ if $G(t,y)$ satisfies the Carath{\'e}odory conditions in $D \ni (0,0)$, i.e. $G$ is measurable in $t$ for each fixed $y$, continuous in $y$ for each fixed $t$ and there is an integrable function $m(t)$ such that
        $$
            \left| G(t,y) \right| \leq m(t)\ \ \text{for any $(t,y) \in D$.}
        $$
        It is easy to see that if the condition of Theorem \ref{p7t2} holds, then $G(x,y) = - F'_x(x,y)/F'_y(x,y)$ satisfies the Carath{\'e}odory conditions in some neighborhood of $(0,0)$. We cannot use the Carath{\'e}odory theorem to provide a simple and direct proof of Theorem \ref{p7t2}.
    \end{remark}

    \subsection{Smoothness of solutions}
    To apply the implicit function theorem, it is necessary to obtain some properties of solutions of the Dirac equation. Recall that $\p(x,\l,t)$ is the fundamental solution of the shifted Dirac equation (\ref{p2e14}) and $\p(x,\l) = \p(x,\l,0)$ is the fundamental solution of the Dirac equation (\ref{p2e1}).

    We say that a matrix-valued function is from $\cH^1(I)$ or $\mH(I_1,I_2)$ if each of its components is from $\cH^1(I)$ or $\mH(I_1,I_2)$. Recall that the dot denotes the derivative with respect to $t$, i.e. $\dot{u} = du/dt$. Let $[A,B] = AB - BA$ be the commutator of two matrices.

    \begin{lemma} \label{p6l1}
        Let $I_1,I_2 \ss \R$ be open bounded intervals and let $F:I_1 \ts I_2 \to \R$ be defined by $F(t,\l) = \p(1,\l, t)$, $(t,\l) \in I_1 \ts I_2$. Then we have $F,F'_{\l} \in \mH(I_1,I_2)$. Moreover, for almost all $t \in \R$ and for each $\l \in \C$ we get
        \begin{align}
            \dot{\p}(1, \l, t) &= \left[ J(V(t)-\l), \p(1, \l, t) \right], \label{p7e5} \\
            \partial_{t} \p'_{\l}(1, \l, t) &= \left[ J(V(t)-\l), \p'_{\l}(1, \l, t) \right] - \left[ J, \p(1, \l, t) \right]. \label{p6e11}
        \end{align}
    \end{lemma}
    \begin{proof}
        We show that $F \in \mH(I_1,I_2)$. By Theorem \ref{p2t1}, each element of $\p(1,\cdot,t)$, $t \in [0,1]$, is an entire function. Since $\p(1,\cdot,t) = \p(1,\cdot,t+1)$, it follows that this statement holds true for any fixed $t \in \R$. Thus, $F,F'_{\l}$ are continuously differentiable function of $\l$ for any fixed $t$.

Using (\ref{p2e4}) and the representation of the inverse matrix
$a^{-1} = \adj a / \det a$, we can rewrite (\ref{p2e15}) as follows
$\p(1,\l,t) =  \p(1+t,\l) \adj \p(t,\l)$. Thus, by Theorem
\ref{p2t1}, each element of $\p(1,\l,\cdot)$ is a linear combination
of functions from $\cH^1(\T)$. It follows that $\p(1,\l,\cdot) \in
\cH^1(I_1)$ for any fixed $\l \in I_2$. Differentiating
(\ref{p2e15}) by $t$ and using (\ref{p2e1}), $V(t+1) = V(t)$, we get
(\ref{p7e5}). Using estimate (\ref{p2e8}) and boundedness of $I_2$,
we get for any $(t,\l) \in I_1 \ts I_2$
        $$
            |\dot{\p}(1,\l,t)| = \left|\left[ J(V(t)-\l), \p(1,\l,t) \right]\right| \leq C_1 \left( \left|q_1(t)\right| + \left|q_2(t)\right| \right) + C_2 = g(t),
        $$
        for some $C_1,C_2>0$. Since $q_1, q_2 \in L^2(\T)$ and $I_1$ is bounded, it follows that $g \in L^2(I_1)$.

        Now we prove that $F'_{\l} \in \mH(I_1,I_2)$. As in the previous case, it is easy to see that $\partial_{\l} \p(1,\cdot,t)$ is entire function for each $t \in I_1$. Differentiating (\ref{p2e15}) by $\l$, we get
        \[ \label{p6e2}
            \begin{aligned}
                \partial_{\l} \p(1,\l,t) &= \p(1+t,\l) \left( -\p^{-1}(t,\l) \left( \partial_{\l} \p(t,\l) \right) \p^{-1}(t,\l) \right)\\
                &+ \left( \partial_{\l} \p(1+t,\l) \right) \p^{-1}(t,\l).
            \end{aligned}
        \]
        Let us prove that $\partial_{\l} \p(\cdot, \l) \in \cH^1(I_1)$. The function $\p(x, \l)$ is a solution of the equation
        $$
            J\p'(x,\l) + V(x) \p(x,\l) = \l \p(x,\l).
        $$
        Differentiating this equation by $\l$ and changing the order of differentiation, we get
        \[ \label{p6e9}
            J \left( \partial_{\l} \p(x,\l) \right)' + V(x) \partial_{\l} \p(x,\l) = \l \partial_{\l} \p(x,\l) + \p(x,\l).
        \]
        It is easy to see that we can change the order of differentiation for smooth potentials. This is possible for potentials from $\cP$, since the left-hand side and the right-hand side of the equation are continuous functions of $V$ and the set of smooth potentials is dense in $\cP$. One can prove that $\p$, and $\partial_{\l} \p$ are continuous functions at $V$ as in \cite{PT} for the Schr{\"o}dinger operator. The solution of (\ref{p6e9}) has the form
        \[ \label{p6e5}
            \partial_{\l} \p(x, \l) = \int_0^x \p(x-\t,\l,\t) J \p(\t,\l) d\t = \p(x,\l) \int_0^x \p^{-1}(\t,\l) J \p(\t,\l) d\t.
        \]
        It is easy to see that the Lebesgue integral of function from $\cH^1(I_1)$ is also function from $\cH^1(I_1)$ if $I$ is bounded interval. It follows that $\partial_{\l} \p(\cdot, \l) \in \cH^1(I_1)$ for any fixed $\l$. Using (\ref{p6e2}), we get that $\partial_{\l} \p(1,\l,\cdot) \in \cH^1(I_1)$ for any fixed $\l$.

        Let us show boundedness of the derivative. Differentiating (\ref{p6e2}) by $t$, we get (\ref{p6e11}). As in the previous case, it is easy to see that there exists $g \in L^2(I_1)$ such that $\left| \partial_{t} \p'_{\l}(1, \l, t) \right| \leq g(t)$ for any $(t,\l) \in I_1 \ts I_2$.
    \end{proof}
    Substituting the potential $V = q_1J_1 + q_2J_2$ in (\ref{p7e5}), we obtain the following explicit formula for the derivative
    \[ \label{p6e1}
        \dot{\p}(1,\l, t) =
        \ma
            \dot{\vt}_1 & \dot{\vp}_1 \\
            \dot{\vt}_2 & \dot{\vp}_2 \\
        \am =
        \ma
            \vp_1(q_1 - \l) - \vt_2(q_1+\l) & 2\vp_1q_2 - 2a\left(q_1 + \l\right) \\
            - 2\vt_2q_2 + 2a \left( q_1 - \l \right) & -\vp_1(q_1 - \l) + \vt_2(q_1 + \l)
        \am,
    \]
    where $a = a(\l,t)$, $q_i = q_i(t)$, $\vt_i = \vt_i(1,\l,t)$, and $\vp_i = \vp_i(1,\l,t)$, $i =1,2$.

    \subsection{Zeros of solutions}
    To analyse zeros of solutions of the Dirac equation (\ref{p2e1}) it is convenient to introduce the Pr{\"u}fer transformation (see e.g. \cite{W87}). Let $y(x,\l)$ be a vector-valued solution of equation (\ref{p2e1}), then we define $\rho(x,\l)$ and $\beta(x,\l)$ by
    \[ \label{p8e1}
        y(x,\l) = \rho(x,\l) \ma \sin \beta (x,\l) \\ \cos \beta (x, \l) \am, \ \ (x, \l) \in \R \ts \C,
    \]
    where $\rho(x,\l) > 0$, $\beta(0,\l) \in [0,\pi)$, and $\beta(\cdot,\l)$ is a continuous function. From (\ref{p8e1}), it follows that $y_1(x_0,\l) = 0$ if and only if $\beta(x_0,\l) = \pi k$, where $k \in \Z$. Substituting (\ref{p8e1}) in equation (\ref{p2e1}) and combining the equations, we obtain the following equation for $\beta(x,\l)$
    \[ \label{p8e2}
        \beta'(x,\l) = q_2(x) \sin (2 \beta(x,\l)) - q_1(x) \cos (2 \beta(x,\l)) - \l.
    \]
    The following lemma describes the behavior of the Pr{\"u}fer angle $\beta(x,\l)$ in a neighborhood of the point $x_0 \in [0,1)$ such that $y(x_0,\l)=0$. We introduce the following function
    $$
        S(x_0,x) = \left|x - x_0 \right|^{1/2} \left| \int_{x_0}^x \left| V(t) \right|^2 dt \right|^{1/2},\ \ x_0,x \in [0,1).
    $$
    Using the H{\"o}lder inequality, we get $\left| \int_{x_0}^x \left|q_i(t) \right| dt \right| \leq S(x_0,x)$ for any $x_0,x \in [0,1)$ and $i = 1,2$. We also introduce the following notation $f(t) = q_1(t) + \l$;
    \begin{lemma} \label{p8l1}
        Let $y$ be a solution of equation (\ref{p2e1}) for some $\l \neq 0$ and let $y_1(x_0,\l) = 0$ for some $x_0 \in [0,1)$. Then $\beta(\cdot,\l)$ is an absolutely continuous function in some neighborhood of $x_0$ and the following asymptotic holds:
        $$
            \beta(x,\l) = \beta(x_0,\l) - \int_{x_0}^x f(t) dt + O\left( S(x_0,x) \int_{x_0}^x \left|f(t) \right|dt \right)\ \ \text{as $x \to x_0$.}
        $$
        Moreover, if $\sign f(x) = \const \neq 0$ for almost all $x$ in some neighborhood of $x_0$, then the following asymptotic holds:
        \[ \label{p8e5}
            \beta(x,\l) = \beta(x_0,\l) - \int_{x_0}^x f(t) dt \left( 1 + O\left(S(x_0,x)\right) \right)\ \ \text{as $x \to x_0$.}
        \]
    \end{lemma}
    \begin{proof}
        If $y(x_0,\l) = 0$, then we have $\beta(x_0,\l) = \pi k$ for some $k \in \Z$. Using (\ref{p8e1}), we get the following formula:
        $$
            \beta(x,\l) = \arcsin \frac{y_1(x,\l)}{y_1(x,\l)^2 + y_2(x,\l)^2}.
        $$
        The functions $y_1(x,\l)$ and $y_2(x,\l)$ are absolutely continuous and $y_1(x,\l)^2 + y_2(x,\l)^2 \neq 0$ for all $x \in [0,1]$. Moreover, $y_1(x_0,\l) = 0$ and $\arcsin x$ is a Lipschitz continuous function in a neighborhood of the point $x = 0$. Since the composition of an absolutely continuous function and a Lipschitz continuous function is absolutely continuous, it follows that $\beta(\cdot,\l)$ is absolutely continuous in some neighborhood of $x_0$. Integrating (\ref{p8e2}), we get
        \[ \label{p6e7}
            \beta(x,\l) - \beta(x_0,\l) = - \int_{x_0}^x f(t) dt + \eta(x_0,x),
        \]
        where
        $$
            \eta(x_0,x) = \int_{x_0}^x q_2(t) \sin (2 \beta(t,\l)) dt + \int_{x_0}^x q_1(t)(1 - \cos (2 \beta(t,\l))) dt.
        $$
        Let us estimate $\eta(x,x_0)$. We introduce the following function
        $$
            M(x_0,x) = \max_{t \in [x_0,x]} \left| \beta(t,\l) - \beta(x_0,\l) \right|,\ \ x_0,x \in [0,1).
        $$
        Let us estimate the maximum of the left-hand side of (\ref{p6e7}). In the proof we consider all asymptotics when $t \to t_0$. Then using the H{\"o}lder inequality and the Taylor expansion of $\sin x$ and $\cos x$ at zero, we have
        $$
            M(x_0,x) = O\left(\left| \int_{x_0}^x \left| f(t) \right| dt \right| + S(x_0,x) M(x_0,x) \left( 1 + M(x_0,x) \right)\right),
        $$
        which yields
        \[ \label{p6e12}
            M(x_0,x) = O\left( \int_{x_0}^x \left| f(t) \right| dt \right).
        \]
        Estimating $\eta(x_0,x)$ as above and using (\ref{p6e12}), we get
        $$
            \eta(x_0,x) = O\left(S(x_0,x) \int_{x_0}^x \left| f(t) \right| dt \right).
        $$
        If $\sign f(x) = \const \neq 0$ for almost all $x$ in some neighborhood of $x_0$, then for such $x$ we have $\left| \int_{x_0}^x \left| f(t) \right| dt \right| = \left| \int_{x_0}^x f(t) dt \right|$, which yields (\ref{p8e5}).
    \end{proof}
    We introduce the following function. Let $y(x,\l)$ be solution of (\ref{p2e1}) for some $\l \in \C$. Then $N_y(\l)$ denote the number of zeros of the function $y_1(\cdot,\l)$ on the interval $[0,1)$. Let also $\rho_y$ and $\beta_y$ denote the Pr{\"u}fer transformation of $y$.

    \begin{lemma} \label{p8l2}
        If for some $n \in \Z$ and for almost all $t \in \T$
        \[ \label{p8e4}
            \sign(q_1(t) + \a_n^-) = \sign(q_1(t) + \a_n^+) = \const \neq 0,
        \]
        then we have
        $$
            N_{\vp}(\l) =
            \begin{cases}
                \left|n\right| &\text{if $\l \in [\a_n^-,\m_n]$, $n > 0$ or $\l \in [\m_n,\a_n^+]$, $n < 0$},\\
                \left|n\right|+1 &\text{if $\l \in (\m_n,\a_n^+]$, $n > 0$ or $\l \in [\a_n^-,\m_n)$, $n < 0$}.
            \end{cases}
        $$
    \end{lemma}
    \begin{proof}
        Let $\vp_1(x_0,\l) = 0$ for some $x_0 \in [0,1)$. From Lemma \ref{p8l1} and condition (\ref{p8e4}), it follows that $\beta_{\vp}(\cdot,\l)$ is a strictly monotone function in a neighbourhood of $x_0$, for each $\l \in [\a_n^-,\a_n^+]$. And it is either strictly increasing or strictly decreasing in a neighborhood of each zero. Then using the strict monotonicity in a neighborhood of each zero and continuity of $\beta_{\vp}(\cdot,\l)$,  we obtain the following formula:
        $$
            N_{\vp}(\l) = \left| \left \lfloor \frac{\beta_{\vp}(1,\l)}{\pi} \right \rfloor \right|,
        $$
        where $\left \lfloor x \right \rfloor = \max \{ i \in \Z \mid i \leq x \}$. Using the formula (3.7) from \cite{Tes}, we get
        $$
            \partial_{\l} \beta_{\vp}(1,\l) = -\frac{\int_0^1 \vp(\t,\l)^2 d\t}{\rho_{\vp}(1,\l)^2} < 0.
        $$
        From this it follows that $N_{\vp}(\l)$ is a constant function on the intervals $[\a_n^-,\m_n)$, and $(\m_n,\a_n^+]$ and it increases by one when $\l$ runs through $\m_n$ in the direction where $\left|\l\right|$ increased.

        Let us show that $\vp_1(\cdot,\m_n)$ has exactly $\left|n\right|$ zeros on the interval $[0,1)$ if the conditions of the theorem hold true. Consider $\vp(\cdot,\m_n)$ as a curve in $\R^2 \sm \{ (0,0) \}$. Using proposition 7.3 from \cite{GreKap09}, and Lemma \ref{p7l1}, Lemma \ref{p6l1}, one can show that this curve and the map $x \mapsto \ma -\sin \pi n x \\ \cos \pi n x \am$ are homotopic. From this it follows that $\beta_{\vp}(1,\m_n) = \pi n$ and $N_{\vp}(\m_n) = \left|n\right|$.
    \end{proof}

    \begin{lemma} \label{p8l3}
        If $\m_n \neq \a_n^{\pm}$ for some $n \in \Z$ and for almost all $t \in [0,1)$
        $$
            \sign(q_1(t) + \a_n^-) = \sign(q_1(t) + \a_n^+) = \const \neq 0,
        $$
        then we have $N_{\p^{\pm}}(\a_n^{\pm}) = \left|n\right|$ for any choice of sign $\pm$.
    \end{lemma}
    \begin{proof}
        If $\m_n \neq \a_n^{\pm}$, then we get $W(\p^{\pm}(x,\a_n^{\pm}),\vp(x,\a_n^{\pm})) = 1$. Using the Pr{\"u}fer transformation (\ref{p8e1}) for $\p^{\pm}$ and $\vp$, we rewrite the Wronskian of these solutions as follows
        $$
            W(\p^{\pm}(x,\a_n^{\pm}),\vp(x,\a_n^{\pm})) = \rho_{\p^{\pm}}(x,\a_n^{\pm}) \rho_{\vp}(x,\a_n^{\pm}) \sin (\beta_{\p^{\pm}}(x,\a_n^{\pm}) - \beta_{\vp}(x,\a_n^{\pm})).
        $$
        Since Wronskian is constant and does not equal zero, and $\beta_{\p^{\pm}}(\cdot,\a_n^{\pm})$, $\beta_{\vp}(\cdot,\a_n^{\pm}) \in C[0,1]$, we get
        \[ \label{p8e3}
            \beta_{\p^{\pm}}(x,\a_n^{\pm}) - \beta_{\vp}(x,\a_n^{\pm}) \in (0,\pi),\ \ x \in [0,1].
        \]
        As in Lemma \ref{p8l2}, we can get the following formula:
        $$
            N_{\p^{\pm}}(\a_n^{\pm}) = \left| \left \lfloor \frac{\beta_{\p^{\pm}}(1,\a_n^{\pm})}{\pi} \right \rfloor - \left \lfloor \frac{\beta_{\p^{\pm}}(0,\a_n^{\pm})}{\pi} \right \rfloor \right|,
        $$
        where $\left \lfloor x \right \rfloor = \max \{ i \in \Z | i \leq x \}$. Using Lemma \ref{p8l2} and identity (\ref{p8e3}), we get
        \[ \label{p6e3}
            N_{\p^{\pm}}(\a_n^-) \in \{\left|n\right| - 1, \left|n\right| \},\ \ N_{\p^{\pm}}(\a_n^+) \in \{\left|n\right|, \left|n\right| + 1 \}.
        \]

        In Theorem \ref{p0t2}, it is shown that if the conditions of the theorem are satisfied, then $\m_n(t)$ runs monotonically around the gap $\g_n^c$, changing sheets when it hits $\a_n^{\pm}$, and $\m_n(0) = \m_n(1)$. It now follows that the number of points such that $\m_n(\cdot)$ equals $\a_n^-$ equals the number of points such that $\m_n(\cdot)$ equals $\a_n^+$. From iv), Theorem \ref{p0t1}, it follows that the number of zeros of $\p^{\pm}_1(\cdot, \a_n^{\pm})$ coincides with the number of points $t_0 \in [0,1)$, where $\m_n(t_0) = \a_n^{\pm}$. Using (\ref{p6e3}), we get $N_{\p^{\pm}}(\a_n^{\pm}) = \left|n\right|$.
    \end{proof}
\footnotesize

\no {\bf Acknowledgments.} Our study was supported by the RSF grant No 18-11-00032.

\end{document}